\documentclass[12pt,a4paper]{amsart}

\usepackage{amsmath}
\usepackage{amssymb}

\usepackage[all]{xy}

\newtheorem{theo}{Theorem}[section]
\newtheorem{thm}[theo]{Theorem}
\newtheorem{lemma}[theo]{Lemma}
\newtheorem{prop}[theo]{Proposition}
\newtheorem{cor}[theo]{Corollary}
\theoremstyle{remark}

\newtheorem{example}[theo]{Example}
\newtheorem{examples}[theo]{Examples}
\newtheorem{remark}[theo]{Remark}

\def\N{{\mathbb N}}
\def\L{{\mathcal L}}

\def\T{{\mathbb T}}
\def\R{{\mathbb R}}
\def\C{{\mathbb C}}

\def\CL{{\mathcal L}}

\def\F{{\mathcal F}}

\def\Z{{\mathbb Z}}

\newcommand{\clsp}{\overline{\operatorname{span}}}

\newcommand{\End}{\operatorname{End}}
\newcommand{\id}{\operatorname{id}}
\newcommand{\range}{\operatorname{range}}
\numberwithin{equation}{section}

\begin{document}

\title[Multi-resolution analyses for Hilbert modules]{Projective multi-resolution analyses arising from\\direct limits of Hilbert modules}
\author[Nadia S. Larsen]{Nadia S. Larsen}
\address{Department of Mathematics, University of Oslo,
P.O. Box 1053, Blindern, N-0316 Oslo, Norway.}
\email{nadiasl@math.uio.no}
\author[Iain Raeburn]{Iain Raeburn}
\address{School of Mathematical and Physical Sciences, University of Newcastle, 
NSW 2308, Australia}
\email{iain.raeburn@newcastle.edu.au}
\date{January 3, 2007}

\begin{abstract}
The authors have recently shown how direct limits of Hilbert spaces can be used to construct multi-resolution analyses and wavelets in $L^2(\R)$. Here they investigate similar constructions in the context of Hilbert modules over $C^*$-algebras. For modules over $C(\T^n)$, the results shed light on work of Packer and Rieffel on projective multi-resolution analyses for specific Hilbert $C(\T^n)$-modules of functions on $\R^n$. There are also new applications to modules over $C(C)$ when $C$ is the infinite path space of a directed graph. 

\end{abstract}

\thanks{This research was supported by the Australian Research Council, the Research Council of Norway, and  
the Faculty of Mathematics and Natural Sciences at the University of Oslo.}

\maketitle

\centerline{Dedicated to the memory of Gert K. Pedersen}

\section*{Introduction}

A multi-resolution analysis for $L^2(\R^n)$ consists of a two-sided sequence of subspaces $V_k$ which are the dilates of a 
single subspace $V_0$ with an orthonormal basis consisting of the integer translates of a single scaling function $\phi$. The Fourier transform of $V_0$ then has an orthonormal basis of the form 
\[
\{e^{2\pi im\cdot x}\phi(x):m\in\Z^n\},
\]
and is the closure of the set $\{f(e^{2\pi ix})\phi(x):f\in C(\T^n)\}$; in other words, $V_0$ is the closure of the free module over $C(\T^n)$ generated by $\phi$. Packer and Rieffel have shown that one can also obtain \emph{projective multi-resolution analyses} in which the 
initial module $V_0$ is the closure of a finitely generated projective module over $C(\T^n)$ \cite{PR2}.

Over the past few years, various authors have realised that Hilbert modules over $C^*$-algebras provide a fertile environment 
for studying multi-resolution analyses, wavelets and frames \cite{fra-lar,rt,wood,PR1,dr}, and Packer and Rieffel worked in that 
context throughout. They considered a specific Hilbert $C(\T^n)$-module $\Xi$ of functions on $\R^n$, and their projective  
multi-resolution analyses consist of the dilations $V_k$ of one Hilbert submodule $V_0$  (see \cite[Definition~4]{PR2}). The spaces $V_k$ are themselves Hilbert modules over $C^*$-algebras of functions on different compact quotients of $\R^n$, and are also projective $C(\T^2)$-modules, though this is not emphasised in \cite{PR2}.

In our previous paper \cite{lr}, we showed that direct limits of Hilbert spaces provide a useful framework for 
Mallat's famous construction of wavelets from a mirror filter \cite{m}. Here we use direct limits of Hilbert modules over 
$C^*$-algebras to produce projective multi-resolution analyses in Hilbert modules. We believe that our methods shed considerable 
light on the constructions in \cite{PR2}, and yield interesting new information about Packer and Rieffel's module $\Xi$. Our 
methods require that we work in the category of Hilbert modules over a fixed $C^*$-algebra, and the spaces in our multi-resolution 
analyses are Hilbert modules in the same category which can be concretely realised using the tensor powers of a single module. When the 
initial module $V_0$ and its complement $W_0$ in $V_1$ are free, we can use our multi-resolution analyses to find orthonormal module bases; when $W_0$ is not free, 
we obtain module frames in the sense of Frank and Larson \cite{fra-lar}. We can in particular write down specific module bases for 
the module $\Xi$. However, many of our constructions are quite general, and we also describe new examples based on the path spaces of finite directed graphs.

A projective multi-resolution analysis for a Hilbert module $X$ over a $C^*$-algebra $A$ is, loosely speaking, an 
increasing sequence of Hilbert submodules $\{V_k\}$ which yields a direct-sum decomposition $X=V_0\oplus\big(\bigoplus_{k=0}^\infty W_k\big)$; we call it projective because one hopes that $V_0$ and the $W_k$ are finitely generated and projective (as they are in \cite{PR2}). In its full generality, our construction starts with a fixed Hilbert 
module $Y$, a correspondence $M$ over $A$, and an isometry $T$ of $Y$ into the balanced tensor product $Y\otimes_A M$. We build a 
direct system with modules $Y\otimes_A M^{\otimes k}$  and maps $T_k:=T\otimes\id$, and then the direct limit module $Y_\infty$ 
comes with a canonical projective multi-resolution analysis in which $V_0$ is a copy of $Y$, $V_k$ is isomorphic to $Y\otimes_A M^{\otimes k}$,  and $W_k$ is isomorphic to the complement of the range of $T_k$ in $Y\otimes_A M^{\otimes (k+1)}$ (see Proposition~\ref{moddecomp}).

When $A$ has an identity $1_A$ and we take $Y$ to be the free module $A_A$, $A\otimes_A M$ is naturally isomorphic to $M$, and the isometries $T:A\to M$ have the form $a\mapsto m\cdot a$ for vectors $m\in M$ such that $\langle m,m\rangle=1_A$; by parallel with the classical case, we then say that $m$ is a \emph{filter}. So each filter $m$ yields a direct limit $M_\infty$ with a projective multi-resolution analysis. For this construction to be useful, we need to be able to identify interesting modules $X$ (such as the module $\Xi$ from \cite{PR2}) as having the form $M_\infty$. As in \cite{lr}, we do this using a dilation operator and a scaling function. Here, the appropriate notion of dilation operator is a Hilbert module isomorphism $D:X\to X\otimes_A M$, and the scaling function is an element $\phi$ of $X$ such that $D\phi=\phi\otimes m$ (see Corollary~\ref{def_of_MRA}).

For modules such as $X=\Xi$, we usually expect the dilation to be an operator from $X$ to itself, and to achieve this we need to choose appropriate correspondences $M$.  We use correspondences $M_L$ which were introduced by Exel in his study of irreversible dynamics \cite{exel, ev}, and which are associated to an endomorphism $\alpha$ of a $C^*$-algebra $A$ and a transfer operator $L$ for $\alpha$. When $A$ is $C(\T)$ or $C(\T^n)$, there are natural transfer operators such that the filters in $M_L$ are the filters arising in wavelet theory and signal processing, and for these filters we can identify the underlying vector spaces of $M_L$ with $A$ and $X\otimes_A M_L$ with $X$. For the classical case in which $A=C(\T)$, $\alpha(f)(z)=f(z^N)$, and $m_0$ is a low-pass filter, the usual dilation operator and scaling function $\phi\in \Xi\subset L^2(\R)$ satisfy our needs, and our construction converts the tensor-product based multi-resolution analysis for $(M_L)_\infty$ into one for $\Xi$ (Example~\ref{classicex}). Our tensor product construction, however, allows us to build an orthonormal basis for the module $\Xi$ (see Example~\ref{basisforXi}).

We have organised our work as follows. In Section~\ref{isometries}, after describing some 
subtleties associated with isometries on Hilbert modules, we define the direct limit of a system of Hilbert modules. Our direct limits are a little different from those in \cite{aa}: it is important for us that the limit has a universal property for maps which are not inner-product preserving. In Section~\ref{MRA_V},  we prove our main result about the existence of projective 
multi-resolution analyses based on a direct limit of Hilbert modules (Theorem~\ref{limit_from_general_V}), and then specialise 
to modules over a unital $C^*$-algebra and systems in which the initial module is $A_A$. 

Next we specialise to the correspondences $M_L$ associated to the endomorphisms of $A=C(C)$ induced by local homeomorphisms of a compact space $C$ (Section~\ref{section_Exel}).
In Section~\ref{frames}, we describe a general procedure for building frames, gradually specialising until we obtain a specific module basis for the module $\Xi$. In Section~\ref{projPR}, we apply the general construction of Section~\ref{MRA_V} with $Y_0$ the non-free projective $C(\T^2)$-module considered by Packer and Rieffel in \cite[\S 4-5]{PR2}. The final result is the same as theirs, but we believe our approach is more systematic, and helps to explain why some of the choices they made are natural. In our final Section~\ref{graphex}, we present a general procedure for constructing Hilbert modules from inverse limits of 
compact spaces, and we use the correspondence associated in \cite{bro} to a finite 
directed graph $E$ to produce a multi-resolution analysis for a module of functions on the two-sided path space of the graph. This example provides further evidence that, in practice, our abstract construction becomes quite concrete.

\subsection*{Conventions} We will be working in the category of right-Hilbert modules over a fixed $C^*$-algebra $A$. If $X$ is such a module, we denote the inner product by $(x,y)\mapsto \langle x,y\rangle$ and the module action by $(x,a)\mapsto x\cdot a$. When we change the module action or inner product, we try to add a subscript to remind ourselves we have done this; for example, if $L$ is a transfer operator for an endomorphism of $A$, we write $\langle \cdot,\cdot\rangle_L$ for the inner product on $X$ defined by $\langle x,y\rangle_L=L(\langle x,y\rangle)$. Some of our Hilbert $A$-modules are correspondences over $A$, which means we also have a left action $(a,x)\mapsto a\cdot x$ given by a homomorphism of $A$ into the $C^*$-algebra $\L(X)$ of adjointable operators on $X$. 

All tensor products of Hilbert modules in this paper are internal tensor products. Thus when we form $X\otimes_A M$, for example, we are assuming that $M$ is a correspondence, and that we have completed the algebraic tensor product $X\odot M$ in the inner product defined on elementary tensors by
\[
\langle x\otimes m,y\otimes n\rangle:=\langle\langle y,x\rangle\cdot m,n\rangle.
\]
Since completing includes modding out by vectors of length zero, it \emph{balances} the tensor product by making $(x\cdot a)\otimes m=x\otimes (a\cdot m)$, and we write $X\otimes_A M$ to remind us of this (see \cite[page~48]{RW}, for example). However, since every tensor product here is balanced, we simplify notation by writing $x\otimes y$ rather than $x\otimes_A y$ for the image of an elementary tensor in $X\otimes_A M$.

Several of our applications involve a specific module $\Xi$ over $C(\T^n)$ considered by Packer and Rieffel in
\cite[\S1]{PR2}. It is defined for any $n\geq 1$, though we are primarily interested in $n=1$ and $n=2$. As a set, $\Xi$ consists of the continuous functions $\xi:\R^n\to \C$ for which there 
is a constant  $K$ such that $\sum_{k\in \Z^n}\vert \xi(t-k)\vert^2\leq K$ for all $t$ in 
$\R^n$, 
and such that the function defined by this sum is continuous; with right action and 
inner product defined by
\begin{equation}
(\xi\cdot f)(t):=\xi(t)f(e^{2\pi i t})
\label{right_action_Xi}
\end{equation}
and 
\begin{equation}
\langle \xi,\eta \rangle(e^{2\pi i t})=\sum_{k\in \Z^n}\overline{\xi(t-k)}\eta(t-k),
\label{def_ip_Xi}
\end{equation}
$\Xi$ is a Hilbert $C(\T^n)$-module \cite[Proposition~7]{PR2}. We freely use key properties of $\Xi$ established by Packer and Rieffel, and especially Propositions~13 and~14 of \cite{PR2}.

\section{Isometries and direct limits of  Hilbert modules}\label{isometries}

For a bounded operator $T:H\to K$ between Hilbert spaces, the following statements are equivalent:
\begin{itemize}
\item[(i)] $T$ is isometric: $\|Th\|=\|h\|$ for every $h\in H$;
\smallskip
\item[(ii)] $T$ preserves the inner product: $(Tg\,|\,Th)=(g\,|\,h)$ for every $g,h\in H$;
\smallskip
\item[(iii)] $T$ satisfies $T^*T=1$.
\end{itemize}
In view of (i), we call an operator with these properties an isometry. 

When $X$ and $Y$ are Hilbert modules over a $C^*$-algebra $A$ and $T:X\to Y$ is a bounded $A$-linear map, 
we easily obtain (iii) $\Longrightarrow$ (ii) $\Longrightarrow$ (i). 
The converse (i) $\Longrightarrow$ (ii) is also true, but is nontrivial: to see this, observe that 
the range of an $A$-module homomorphism $T$ satisfying (i) is automatically a closed submodule of 
$Y$, and \cite[Theorem~3.5]{lance} implies that it is unitary as a map into $T(X)$, hence 
inner-product preserving. On the other hand, it is no longer true that (ii) $\Longrightarrow$ (iii), unless we know that $T$ is adjointable (see Lemma~\ref{Hilbert_module_isometry} below). Not all inner-product preserving maps are adjointable: 
the standard counterexample is the inclusion of $C_0((0,1])$ in $C([0,1])$, viewed as a map of 
Hilbert $C([0,1])$-modules (see \cite[page~21]{lance} or \cite[Example~2.19]{RW}). We are 
interested in both adjointable and non-adjointable inner-product preserving operators, and for 
us the key difference is that adjointable inner-product preserving maps have complemented range (see Lemma~\ref{Hilbert_module_isometry}).

A closed submodule $M$ of a Hilbert $A$-module $X$ is 
\emph{complemented} if the map $(m, n)\mapsto m+n$ is a Hilbert 
module isomorphism of $M\oplus M^\perp$ onto $X$.  
When this is the case, there is an orthogonal projection $P:X\to M$, 
and $P$ is an adjointable operator on $X$ such that 
$P^2=P=P^*$. The converse observation is also valid: if $P\in \mathcal{L}(X)$ satisfies $P^2=P=P^*$ then 
$M:=P(X)$ is complemented with complement $M^\perp=(1-P)(X)$. The next lemma sums up the 
key properties: it is similar to but not quite the same as \cite[Proposition~3.6]{lance}.

\begin{lemma}\label{Hilbert_module_isometry}
Let $A$ be a $C^*$-algebra, and suppose that $X$ and $Y$ are Hilbert $A$-modules 
and $S:X\to Y$ is inner-product preserving. Then $S$ is an $A$-module homomorphism, 
and $S$ is adjointable if and only if 
the range of $S$ is a complemented $A$-submodule of $Y$. If this is the case,  
then $S^*S=1$, and $SS^*$ is the orthogonal projection onto the range of $S$.
\end{lemma}

\begin{proof}
That $S$ is an $A$-module homomorphism follows by computing $\|S(x\cdot a) -
(Sx)\cdot a\|^2$ for $x\in X$ and $a\in A$. If $S$ is adjointable, then 
$S^*S=1$, and $P:=SS^*$ is a 
self-adjoint projection in $\mathcal{L}(Y)$ with $P(Y)=S(X)$. We deduce from 
the observations in the paragraph preceding the lemma that $S(X)$ is 
complemented. 
Conversely, if $S(X)$ is complemented and $P$ is the projection onto $S(X)$, then the inverse $S^{-1}:S(X)\to X$ satisfies 
\[
\langle Sx, y\rangle=\langle PSx, y\rangle=\langle Sx, Py\rangle
=\langle Sx, SS^{-1}Py\rangle
=\langle x, S^{-1}Py\rangle,
\]
and hence $S^{-1}P$ is an adjoint for $S$.
\end{proof}

In view of the distinction between (ii) and (iii) it is potentially ambiguous to talk about ``isometries of Hilbert bimodules'', and we try to be precise at least in the formulations of our results. As an example: the following easy lemma gives one way to construct isometries.

\begin{lemma}\label{def_of_Sx}
Let $A$ be a $C^*$-algebra with identity and $X$ a Hilbert $A$-module. For every
$x\in X$, the map $S_x:a\to x\cdot a$ of $A_A$ into $X$ is adjointable with adjoint given by 
$S_x^*y=
\langle x,y \rangle$. The map $S_x$ is inner-product preserving if and only if 
$\langle x,x \rangle=1$.
\end{lemma}

We now describe our direct limit and its universal property.

\begin{prop}\label{dir_lim}
Suppose that $A$ is  a  $C^*$-algebra and 
$$
X_0\overset{T_0}{\longrightarrow} X_1\overset{T_1}{\longrightarrow}
X_2 \longrightarrow\dots
$$
is a direct system of Hilbert $A$-modules $X_k$ in which each
$T_k$ is inner-product preserving (but not necessarily adjointable). 
\smallskip 

\textnormal{(a)} There are a Hilbert $A$-module $X_\infty$ and inner-product 
preserving maps $\iota^k:X_k\to X_\infty$ with the following property: whenever 
$R_k:X_k\to Z$ are bounded $A$-module homomorphisms of $X_k$ into another Hilbert 
$A$-module $Z$ such that $R_{k+1}\circ T_{k}=R_{k}$ and $\|R_k\|\leq M$ for all $k\geq 0$, 
there is a unique bounded $A$-module homomorphism 
$R_\infty$ such that $\|R_\infty\|\leq M$ and $R_\infty\circ \iota^k=R_k$ for 
$k\geq 0$. We illustrate this with the commutative diagram:
\[\xygraph{
{X_0}="v0":[rr]{X_{1}}="v1"_{T_0}:[rr]{X_{2}}="v2"_{T_{1}}:[rr]{\cdots}="v3"_{T_2}:@{}[rrr]
{X_\infty}="v4"
"v0":@/^{40pt}/^{\iota^0}"v4""v1":@/^{20pt}/"v4"^{\iota^1}"v2":@/^{7pt}/"v4"^{}
"v1":[dd]{Z}="w1"_{R_1}"v2":"w1"^{R_2}"v0":"w1"_{R_0}"v4":@{-->}"w1"^{R_\infty}
}\]
If every $T_k$ is adjointable, so is every $\iota^k$, and if every 
$R_k$ is inner-product preserving, so is $R_\infty$.
\smallskip

\textnormal{(b)} The pair $(X_\infty,\iota^k)$ is essentially unique: whenever 
$(Y,j^k)$ is a pair with the property described in \textnormal{(a)}, the space 
$\bigcup_{k=1}^\infty j^k(X_k)$ is dense in $Y$, and 
there is an isomorphism $\theta$ of $X_\infty$ onto $Y$ such that $\theta\circ
\iota^k=j^k$ for all $k$.
\end{prop}

The pair $(X_\infty,\iota^k)$ is called the \emph{direct limit} of the system 
$(X_k,T_k)$, and we usually write $X_\infty=\varinjlim X_k$ or $(X_\infty,
\iota^k)=\varinjlim(X_k,T_k)$. We refer to the property in part (a) as the 
\emph{universal property of the direct limit}.

\begin{proof}
To construct such a module, we take the quotient $Q$ of the disjoint union 
$\bigsqcup_{k=0}^\infty X_k$ by the equivalence relation which identifies each $x\in X_k$ 
with $T_kx\in X_{k+1}$. Since the maps $T_k$ are $A$-module homomorphisms and 
preserve the inner products, the quotient $Q$ is an $A$-module and carries an 
$A$-valued inner product. We take $X_\infty$ to be the completion of $Q$ in the 
norm defined by the inner product, which is naturally a Hilbert $A$-module. We write 
$[x]$ for the class of $x\in X_k$ in $Q$ or $X_\infty$, and define $\iota^kx=[x]$ for $x\in X_k$. 

Given $R_k:X_k\to Z$ as in (a), elementary algebra shows that 
there is a well-defined $A$-module homomorphism $R_\infty$ of $Q$ into $Z$ such that 
$R_\infty([x])=R_kx$ for $x\in X_k$. The norm of $R_\infty([x])$ satisfies
\[
\|R_\infty([x])\|=\|R_kx\|\leq \|R_k\|\,\|x\|=\|R_k\|\,\|[x]\|\leq M\|[x]\|,
\]
so $R_\infty$ is bounded on $Q$ and hence extends to a bounded linear operator 
$R_\infty$ with the required properties. Any bounded operator $T:X_\infty\to Z$ satisfying $T\circ\iota^k=R_k$ agrees with $R_\infty$ on $Q$, and hence by continuity also on $X_\infty=\overline Q$. If every $R_k$ is inner-product preserving, then $R_\infty$ is inner-product preserving on $Q$, and the continuity of the inner products implies that $R_\infty$ is inner-product preserving on $X_\infty$.

Now suppose that each $T_k$ is adjointable, and note that $\| T_k^*\|=\|T_k\|=1$ (by, for example, \cite[Proposition 1.10]{ekqr}). To verify the statement about adjointability of $\iota^k$ we apply the universal property to
\[\xygraph{
{X_k}="v0":[rr]{X_{k+1}}="v1"^{T_k}:[rr]{X_{k+2}}="v2"^{T_{k+1}}:[rrr]{\cdots\cdots}="v3"^{T_{k+2}}:[rrr]
{X_\infty}="v4"
"v1":[dd]{X_k}="w1"^{T_k^*}"v2":"w1"^{T_k^*T_{k+1}^*}"v0":"w1"_{\id}"v4":@{-->}"w1"
}\]
to get a bounded $A$-module homomorphism $R_\infty:X_\infty\to X_k$. Then for $x\in X_k$ and $y=
\iota^{k+n}z\in 
X_\infty$, we have
\begin{align*}
\langle x,R_\infty y\rangle
&=\langle x,(R_\infty \circ\iota^{k+n})z\rangle\\
&=\langle x,T_k^*T_{k+1}^*\cdots T_{k+n-1}^*z\rangle\\
&=\langle T_{k+n-1}\cdots T_{k+1}T_kx,z\rangle\\
&=\langle \iota^{k+n}(T_{k+n-1}\cdots T_{k+1}T_kx),\iota^{k+n}z\rangle\\
&=\langle \iota^kx,y\rangle.
\end{align*}
Since both sides are continuous in $y$, this extends to all $y\in X_\infty$, and $R_\infty$ is an adjoint for $\iota^k$, as claimed.

To prove (b), we let $Z:=\overline{\bigcup_{k=1}^\infty j^k(X_k)}$ and aim to prove that $Y=Z$. 
Viewing the maps $j^k$ as having range in $Z$, and applying the universal property to $j^k:X_k\to Z$, gives a 
map $j^\infty:Y\to Z$ such that $j^\infty\circ j^k=j^k$. But now we can view $j^\infty$ as a map of 
$Y$ into $Y$ such that we have a commutative diagram
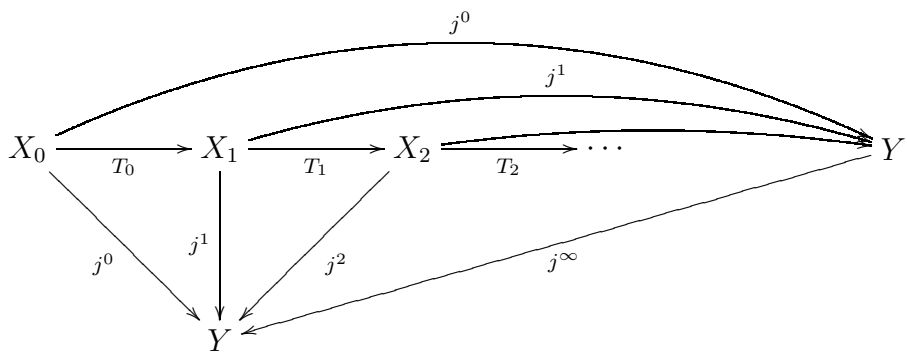
\[\xygraph{
{X_0}="v0":[rr]{X_{1}}="v1"_{T_0}:[rr]{X_{2}}="v2"_{T_{1}}:[rr]{\cdots}="v3"_{T_{2}}:@{}[rrr]
{Y}="v4"
"v0":@/^{40pt}/^{j^0}"v4""v1":@/^{20pt}/"v4"^{j^1}"v2":@/^{7pt}/"v4"^{}
"v1":[dd]{Y}="w1"_{j^1}"v2":"w1"^{j^2}"v0":"w1"_{j^0}"v4":"w1"^{j^\infty}
}\]
Since the diagram also commutes when we replace $j^\infty$ by the identity map $\id_Y$ from $Y$ to $Y$, 
the uniqueness in the universal property forces $j^\infty=\id_Y$ and $Y=Z$.

The homomorphism $\theta$  is obtained by applying the universal property of $(X_\infty,
\iota^k)$ to the maps $j^k:X_k\to Y$. It is an isomorphism because applying the universal property of $(Y,j^k)$ 
to $\iota^k:X_k\to X_\infty$ gives an inverse.
\end{proof}

\begin{remark}
It is tempting to think that we are working exclusively in the category of Hilbert modules and 
inner-product preserving maps. However, in proving that $\iota^k$ is adjointable, we had to apply 
the universal property of Proposition~\ref{dir_lim} to maps which do not preserve the inner product. This means we cannot directly apply the results of \cite{aa}, for example.
\end{remark}

\begin{examples}
For an example where each $T_k$ is inner-product preserving but not adjointable, 
take $T_k$ to be the map of $C_0((-k,k))$ into $C_0((-(k+1),k+1))$ which extends 
$f\in C_0((-k,k))$ to be $0$ outside $(-k,k)$ (see \cite[Example 2.19]{RW}). 
When we view each $C_0((-k,k))$ as a Hilbert $C_0(\R)$-module, the direct limit is isomorphic to $C_0(\R)$, 
with the inclusion maps playing the role of the $\iota^k$. 
When we view each $C_0((-k,k))$ as a Hilbert $C(\T)$-module with 
\[
(f\cdot a)(x)=f(x)a(e^{2\pi ix})\ \text{ and }\ \langle f,g\rangle(e^{2\pi ix})=\sum_{l=-k+1}^{k}
\overline{f(x-l)}g(x-l) \text{ for $x\in [0,1)$,}
\]
then the direct limit is the Hilbert $C(\T)$-module $\Xi$ constructed in \cite[\S1]{PR2}. To see this, 
we recall from \cite[Proposition~3]{PR2} that $\Xi$ is complete, and apply Proposition~\ref{dir_lim} 
with $R_k$ the inclusion of $C_0((-k,k))$ in $\Xi$. The induced map $R_\infty$ has dense range by 
\cite[Proposition~4]{PR2}, and hence is surjective.
\end{examples}

\begin{example} For an example in which every $R_k$ is adjointable but $R_\infty$ is not, let 
$A=C(\N\cup \{\infty\})$, let $X_k=\operatorname{span}\{
\delta_j: 0\leq j \leq k\}\subset A$, and let $T_k:X_k\hookrightarrow X_{k+1}$ 
be the inclusion. Then with $X_\infty=\clsp\{\delta_j: j\geq 0\}=c_0(\N)$ and 
$\iota^k$ the inclusion maps, $(X_\infty,\iota^k)$ has the universal property. By 
uniqueness, the map $R_\infty$ associated to the inclusions $R_k:
X_k\hookrightarrow A$ has to be the inclusion of $c_0(\N)$ in $A$. However, $R_\infty(c_0(\N))^\perp 
=\{0\}$, so the range of $R_\infty$ is not complemented in $A$, and Lemma~\ref{Hilbert_module_isometry} 
implies that $R_\infty$ is not adjointable.
\end{example}

\section{Multi-resolution analyses from direct limits of modules}\label{MRA_V}

The direct systems of interest to us involve the tensor powers of a fixed correspondence $M$ over $A$. We have an important motivating example in mind.

\begin{example}\label{keyexonT}
In this example, the $C^*$-algebra $A$ is $C(\T)$, and the module $M$ also has underlying space $C(\T)$. (We will try to distinguish elements of the coefficient algebra by calling them $a$ or $b$.) The module structure depends on a fixed integer $N\geq 2$; the actions are defined by
\[
(f\cdot a)(z):=f(z)a(z^N)\ \text{ and }\ (a\cdot f)(z)=a(z)f(z)\ \text{ for $f\in M$ and $a\in A$,}
\]
and the inner product is given by 
\[
\langle f,g\rangle(z)=\frac{1}{N}\sum_{\{w\in\T\;:\;w^N=z\}}\overline{f(w)}g(w).
\]
It is easy to check that the operator $f\mapsto a\cdot f$ is then adjointable with adjoint $g\mapsto a^*\cdot g$, 
so the left action gives a $C^*$-algebra homomorphism of $A$ into the algebra $\L(M)$ of adjointable operators on the right Hilbert $A$-module $M$.
\end{example} 

To construct our direct systems, we start with a fixed Hilbert $A$-module $Y$, and an isometry $T$ of $Y$ into the balanced tensor product $Y\otimes_A M$ by our fixed correspondence $M$. Again we have a specific motivating example.

\begin{example}\label{excorresponA}
Suppose that $A$ is unital with identity $1_A$, the correspondence $M$ is essential in the sense that 
$1_A\cdot f=f$ for $f\in M$, and take $Y=A_A$. Recall that a \emph{filter} in $M$ is an element $m$ such that $\langle m,m\rangle=1_A$. Then for every filter $m\in M$, the isometry $S_m:A\to M$ of Lemma~\ref{def_of_Sx} fits our model. (To see this, we 
just need to observe that the map $a\otimes_A n\mapsto a\cdot n$ is an isomorphism of $A\otimes_A M$ onto $M$, so we can view $S_m$ as a map into $A\otimes_A M$.) Filters in our motivating Example~\ref{keyexonT} are functions $m:\T\to\C$ such that
\[
\frac{1}{N}\sum_{w^N=z}|m(w)|^2=1\ \text{ for every $z\in\T$,}
\]
or in other words, filters in the sense of wavelet theory (see \cite[Equations (1.14) and (1.25)]{Bra-Jor-book}, for example),
and then $S_m:A\to M$ is defined by
\[
(S_ma)(z)= m(z)a(z^N).
\]
\end{example}

The other maps $T_k$ in our direct system $(X_k,T_k)$ of Hilbert modules will be formed 
by tensoring $T:Y\to T\otimes_A M$ with the identity maps $\id_k$ on the tensor powers $M^{\otimes k}:=M\otimes_A\cdots\otimes_A M$. There are some subtleties involved in forming such tensor products, so we briefly discuss this construction.

Suppose that $Y$ and $Z$ are Hilbert modules over a $C^*$-algebra $A$, and $M$ is a correspondence over $A$. 
For every $T\in\mathcal{L}(Y,Z)$, there is a unique adjointable map $T\otimes\id$ from the 
internal tensor product $Y\otimes_A M$ to $Z\otimes_A M$ characterised by
\begin{equation}
(T\otimes \id)(y\otimes m)=(Ty)\otimes m
\label{def_Q_tensor_id}\ \text{ for $y\in Y$ and $m\in M$.}
\end{equation}
To prove that there is such a map $T\otimes\id$ requires 
non-trivial arguments 
(see \cite[page~42]{lance}), 
but the characterising property \eqref{def_Q_tensor_id} makes it easy to manipulate. 
For example, one can check by computing on elementary tensors that $(S\circ T)\otimes\id=(S\otimes\id)
\circ(T\otimes\id)$, and that $(T\otimes\id_M)\otimes\id_N=T\otimes\id_{M\otimes_A N}$. 
Adjointability plays a crucial role in proving that $T\otimes\id$ extends to the completion, so 
we cannot in general form $T\otimes\id$  for non-adjointable operators $T$, even if they are
norm-bounded. However, if $T:Y\to Z$ is inner-product preserving, then one can verify directly 
that there is a well-defined linear operator $T\otimes \id$ on the algebraic tensor product $Y\odot M$ which 
preserves the internal tensor-product norm, and hence induces an inner-product preserving map 
$T\otimes\id$ on $Y\otimes_A M$ satisfying \eqref{def_Q_tensor_id}.

Our construction involves a Hilbert-module isomorphism $U$ of $Y$ onto $Y\otimes_A M$, and the 
isomorphisms $U^k$ defined inductively by $U^0=U$ and 
\begin{equation}\label{defDk}
U^{k+1}=(U\otimes {\id_k})\circ U^k:Y\to Y\otimes_A M^{\otimes (k+1)}.  
\end{equation}
We can verify by calculations on elementary tensors that
\begin{equation}\label{altdefDk}
U^{k+1}=(U^{k}\otimes\id_1)\circ U.
\end{equation}

We can now formulate the main result of the section.

\begin{thm}\label{limit_from_general_V}
Suppose that $M$ is a correspondence over a $C^*$-algebra $A$, $Y$ is a Hilbert $A$-module, and $T:Y\to Y\otimes_A M$ is inner-product 
preserving and adjointable. Define $T_k:=T\otimes \id_k:Y\otimes_AM^{\otimes k}\to Y\otimes_AM^{\otimes (k+1)}$.
\smallskip

\textnormal{(a)} Let $(Y_\infty, \iota^k)=\varinjlim (Y\otimes_A M^{\otimes k}, T_k)$. Then 
there is a Hilbert-module isomorphism $R$ of $Y_\infty$ 
onto $Y_\infty\otimes_A M$ such that $R\circ \iota^{k+1}=\iota^k \otimes 
\id$ for $k\geq 0$. The submodules $Y_k:=\iota^k(Y\otimes_A M^{\otimes 
k})$ of $Y_\infty$ satisfy:
\begin{enumerate}
\item $Y_0$ is a complemented $A$-submodule of $Y_\infty$;
\smallskip
\item $Y_k=R^{-k}(Y_0\otimes_A M^{\otimes k})$ for $k\geq 0$;
\smallskip
\item $Y_k$ is a complemented submodule of $Y_{k+1}$ for $k\geq 0$;
\smallskip
\item $\bigcup_{k=0}^\infty Y_k$ is dense in $Y_\infty$.
\end{enumerate}

\textnormal{(b)} Suppose $X$ is a Hilbert $A$-module, $D:X\to X\otimes_A M$ is 
a Hilbert-module isomorphism and $Q_0:Y\to X$ is an inner-product 
preserving map such that $D\circ Q_0=(Q_0\otimes \id)\circ T$. 
Then there is 
an inner-product preserving map $Q:Y_\infty\to X$ such that $Q\circ \iota^0=Q_0$ 
and $D\circ Q=(Q\otimes\id)\circ R$. The submodules $V_k:=Q(Y_k)$ of $Q(Y_\infty)$ satisfy 
\begin{equation}\label{formforV_k}
V_k=D^{-k}\circ(Q_0\otimes\id_k)(Y\otimes_A M^{\otimes k})
\end{equation}
and have properties mirroring those of \textnormal{(i)--(iv)}.
\end{thm}

Since the modules in our analysis are often finitely generated and projective, we call a sequence of submodules $\{Y_k\}$ satisfying (i)--(iv) a \emph{projective multi-resolution analysis}  
for $Y_\infty$. Part (b) then says that the sequence $\{V_k\}$ is a projective multi-resolution analysis for the (necessarily closed) submodule $Q(Y_\infty)$ of $X$. We will explain in Example~\ref{classicex} how these projective resolution analyses give the ones considered by Packer and Rieffel in \cite{PR2}.

\begin{proof}
(a) Since $T$ is inner-product preserving and adjointable, so is each $T_k$, and we can form the 
direct limit $(Y_\infty,\iota^k)$; it follows from Proposition~\ref{dir_lim} that each $\iota^k$ is  
inner-product preserving and adjointable. For $k\geq 1$, $R_{k+1}:=\iota^k\otimes\id$ is an 
inner-product preserving map from 
$Y\otimes_A M^{\otimes (k+1)}$ into $Y_\infty\otimes_A M$. Since 
\begin{align*}
R_{k+1}\circ T_k
&=(\iota^{k}\otimes \id)\circ(T\otimes\id_k)=(\iota^k\otimes\id)\circ(T_{k-1}\otimes \id)\\
&=(\iota^k\circ T_{k-1})\otimes\id=\iota^{k-1}\otimes \id=R_k,
\end{align*}
we have a commutative diagram
\[\xygraph{
{Y}="v0":[rr]{Y\otimes_A M}="v1"_{T}:[rr]{Y\otimes_A M^{\otimes 2}}
="v2"_{T_1}:[rr]{\cdots}="v3"_{T_2}:@{}[rrr]
{Y_\infty}="v4"
"v0":@/^{40pt}/^{\iota^0}"v4""v1":@/^{20pt}/"v4"^{\iota^1}"v2":@/^{7pt}/"v4"^{}
"v1":[dd]{Y_\infty\otimes_A M}="w1"_{R_1}"v2":"w1"^{R_2}"v4":@{-->}"w1"^{R}
}\]
in which every solid arrow is inner-product preserving. Thus the universal property 
of $(Y_\infty, \iota^k)$ gives an inner-product preserving map $R:=R_\infty:Y_\infty\to Y_\infty\otimes_A M$ satisfying 
$R\circ \iota^{k}=\iota^{k-1} \otimes \id$.  
Since elements of the form $\iota^{k-1} b$ as $k$ and $b$ vary are 
dense in 
$Y_\infty$,  $R$ is surjective, and hence an isomorphism of Hilbert modules.

To prove (i), note that $Y_0$ is the range of the inner-product preserving and adjointable map 
$\iota^0$,  and hence by Lemma~\ref{Hilbert_module_isometry} is complemented in $Y_\infty$. 
To prove (ii), we note that it is trivially true for $k=0$; for $k\geq 1$ we use \eqref{altdefDk} to compute
\begin{align*}
R^{k}Y_k&=R^{k}\circ \iota^k(Y\otimes_A M^{\otimes k})\\
&=(R^{k-1}\otimes\id_1)\circ(R\circ \iota^k)(Y\otimes_A M^{\otimes k})\\
&=(R^{k-1}\otimes\id_1)\circ(\iota^{k-1}\otimes \id)(Y\otimes_A M^{\otimes k})\\
&=R^{k-1}(Y_{k-1})\otimes_A M,
\end{align*}
and (ii) then follows from an induction argument. For (iii), we observe that $T_k$ is adjointable, and hence its range is a complemented 
submodule of $Y\otimes_A M^{\otimes (k+1)}$. Since $\iota^{k+1}\circ T_k=\iota^k$, 
$\iota^{k+1}$ maps the range of $T_k$ onto $Y_k=\range\iota^k$, and hence $Y_k$ is complemented 
in $Y_{k+1}$, as claimed. Part (iv) follows from Proposition~\ref{dir_lim}(b).

(b) For every $k\geq 1$, the map 
$Q_k:=D^{-k}\circ(Q_0\otimes \id_k):Y\otimes_A M^{\otimes k}\to X$ is inner-product preserving because $D$ and $Q_0$ are. From 
\eqref{defDk} and the given property of $Q_0$, we deduce that 
\begin{align*}
Q_{k+1}\circ T_k
&=D^{-(k+1)}\circ(Q_0\otimes \id_{k+1})\circ(T\otimes \id_k)\\
&=D^{-k}\circ(D\otimes \id_k)^{-1}\circ((Q_0\otimes \id) \otimes \id_k)
\circ(T\otimes \id_k)\\
&=D^{-k}\circ((D^{-1}\circ(Q_0\otimes \id)\circ T)\otimes \id_k)\\
&=D^{-k}\circ(Q_0\otimes \id_k)\\
&=Q_k.
\end{align*}
Hence the universal property of $(Y_\infty,\iota^k)$ gives an inner-product 
preserving map $Q$ of $Y_\infty$ into $X$ such 
that $Q\circ \iota^k=Q_{k}=D^{-k}\circ(Q_0\otimes \id_k)$. 
Then, using \eqref{altdefDk} again, we have
\begin{align}
(D\circ Q)\circ \iota^{k+1}
&=D\circ Q_{k+1}=D\circ D^{-(k+1)}\circ(Q_0\otimes \id_{k+1})\label{2.4a}\\
&=D\circ D^{-1}\circ(D^{-k}\otimes \id)\circ((Q_0\otimes \id_{k})\otimes \id)\notag\\
&=(D^{-k}\circ(Q_0\otimes\id_k))\otimes \id;\notag
\end{align}
on the other hand, the relation $R\circ \iota^{k+1}=\iota^k\otimes \id$ gives
\begin{align}\label{rhs_Q}
(Q\otimes \id)\circ R\circ \iota^{k+1}=(Q\circ \iota^k)\otimes \id=(D^{-k}\circ(Q_0\otimes \id_k))\otimes \id,
\end{align} 
and \eqref{2.4a} and \eqref{rhs_Q} imply that $D\circ Q=(Q\otimes \id)\circ R$. The properties of $V_k:=Q(Y_k)$ follow 
from those of $Y_k$ because $Q$ is inner-product preserving. An induction argument using the formula \eqref{defDk}
shows that $D^k\circ Q=(Q\otimes \id_k)\circ R^k$, and this implies \eqref{formforV_k}. 
\end{proof}

\begin{prop}\label{moddecomp} With the notation of Theorem~\ref{limit_from_general_V}, let 
$Z_k$ denote the complement~$Y_{k+1}\ominus Y_k$ of $Y_k$ in $Y_{k+1}$. Then there is a natural Hilbert-module direct-sum decomposition
\begin{equation}\label{decompMinfty}
Y_\infty=Y_0\oplus\big(\textstyle{\bigoplus_{k=0}^\infty Z_k}\big),
\end{equation}
and the isomorphism $R^k:Y_\infty\to Y_\infty\otimes_A M^{\otimes k}$ induced by 
$R$ restricts to an isomorphism of $Z_k$ onto $Z_0\otimes_A M^{\otimes k}$. 
\end{prop}
 
For the last part we need a simple lemma:

\begin{lemma}\label{simple_lemma}
Suppose $V_1$ is a complemented submodule of a Hilbert $A$-module $V_2$ and $E$ 
is a correspondence  over $A$. Then $V_1\otimes_A E$ is a complemented submodule of $V_2\otimes_A E$ with 
\[
(V_1\otimes_A E)^\perp=V_1^\perp\otimes_A E.
\] 
\end{lemma}

\begin{proof}
Since the inclusions of $V_1$ and $V_1^\perp$ in $V_2$ induce isometric embeddings of 
$V_1\otimes_A E$ and $V_1^\perp\otimes_A E$ into $V_2\otimes_A E$,  the assertion 
at least makes sense. Let $P\in \CL(V_2)$ be the orthogonal projection of $V_2$ on $V_1$. 
Then since $T\mapsto T\otimes \id$ is a homomorphism of $\CL(V_2)$ into $\CL(V_2\otimes_A E)$,  
$P\otimes\id_E$ is a projection. It is easy to see that the range of $P\otimes\id_E$ contains 
$V_1\otimes_A E$. On the other hand, if $x\in (P\otimes\id_E)(V_2\otimes_A E)$, and $x\sim 
\sum_i v_i\otimes y_i \in V_2\odot E$, then
\[
x=(P\otimes \id_E)x\sim \sum (Pv_i)\otimes y_i\in V_1\otimes_A E,
\]
and $x$ belongs to the closed submodule $V_1\otimes_A E$. Thus $V_1\otimes_A 
E=(P\otimes\id_E)(V_2\otimes_A E)$. Then the projection on $(V_1\otimes_A E)^\perp$ 
is $\id_{V_2\otimes_A E}-P\otimes \id_E=(\id_{V_2}-P)\otimes\id_E$, which by a 
similar argument has range $((\id-P)V_2)\otimes_A E=V_1^\perp\otimes_A E$.
\end{proof}

\begin{proof}[Proof of Proposition~\ref{moddecomp}]
The Hilbert-module direct sum is defined in \cite[Page~6]{lance}. Thus the right hand side of \eqref{decompMinfty} is 
\[
\big\{(y,\{z_k\}_{k\geq 0}): y\in Y_0, z_k\in Z_k, \text{ and }\textstyle{\sum}\langle z_k, 
z_k\rangle \text{ converges in }A\big\},
\]
with inner product given by 
\[
\langle (y, \{z_k\}), (y', \{z'_k\}) \rangle=
\langle y, y'\rangle + \sum_k \langle z_k, z'_k\rangle.
\]
For each $(y,\{z_k\}_{k\geq 0})$ in $Y_0\oplus\big(\textstyle{\bigoplus_{k=0}^\infty 
Z_k}\big)$, the series $\sum\langle z_k, z_k\rangle$ converges in $A$, and hence 
$\sum_k z_k$ converges in $Y_\infty$. Thus  $\Sigma:(y,\{z_k\}_{k\geq 0})
\mapsto y +\sum_k z_k$ is a well-defined map 
from $Y_0\oplus\big(\textstyle{\bigoplus_{k=0}^\infty Z_k}\big)$ into $Y_\infty$. 
This map is inner-product preserving, and its range contains $\bigcup_k Y_k$, and hence is dense in $Y_\infty$; since the range of an inner-product preserving map is closed, $\Sigma$ is surjective. Thus $\Sigma$ is an isomorphism of Hilbert modules, and this is precisely what \eqref{decompMinfty} means.

We know from Theorem~\ref{limit_from_general_V}  that $R$ is an isomorphism of 
$Y_{k+1}$ onto $Y_k\otimes_A M$, and that it carries the submodule $Y_{k}$ onto 
$Y_{k-1}\otimes_A M$. 
Thus it takes $Z_k:=Y_{k+1}\ominus Y_k$ onto the complement of 
$Y_{k-1}\otimes_A M$ in $Y_k\otimes_A M$, which by Lemma~\ref{simple_lemma} 
is $Z_{k-1}\otimes_A M$. An induction 
argument now shows that $R^k=(R^{k-1}\otimes\id)\circ R$ 
carries $Z_k$ onto $Z_0\otimes_A M^{\otimes k}$, as claimed.
\end{proof}

We aim to apply Theorem~\ref{limit_from_general_V} with $Y$ the free module $A_A$ over a unital $C^*$-algebra and $T$ the isometry $S_m$ associated to a filter $m$ in a correspondence $M$ over $A$ (see 
Example~\ref{excorresponA}). 
When $M$ is essential as a left $A$-module, the left actions of $A$ 
give natural isomorphisms of $A\otimes_A M^{\otimes k}$ onto $M^{\otimes k}$, and 
under these isomorphisms the maps $T\otimes \id_k$ from $A\otimes_A 
M^{\otimes k}$ into $(A\otimes_A M)\otimes_A M^{\otimes k}$ become the maps $T_k$ defined by
\begin{equation}
T_kn=T_k(1\otimes n)=(m\cdot 1)\otimes n= m\otimes n\ \text{ for $n\in M^{\otimes k}$.}
\label{new_Sk}
\end{equation}
Let $(M_\infty, \iota^k):=\varinjlim (M^{\otimes k}, 
T_k)$. Then Theorem~\ref{limit_from_general_V}(a) says that there is an 
isomorphism $R:M_\infty\to M_\infty\otimes_A M$ characterised by $R\circ \iota^{k+1}=
\iota^k \otimes\id$ for $k\geq 0$, and that the submodules $Y_k:=\iota^k(M^{\otimes k})$ form a projective multi-resolution analysis for $M_\infty$ with $Y_0\cong A_A$.

To find concrete implementations of the pair $(M_\infty, R)$ using part (b) of 
Theorem~\ref{limit_from_general_V}, we need a Hilbert $A$-module $X$, 
an isomorphism $D:X\to X\otimes_A M$, and an inner-product preserving map 
$Q_0:A\to X$ such that $D\circ Q_0=(Q_0\otimes \id)\circ S_m$. The map $Q_0$ 
is determined by its value $\phi:=Q_01$; notice that $\phi$ must satisfy 
$\langle \phi, \phi\rangle=1$ and, remembering that $m\in M$ identifies with 
$1\otimes m$ in $A\otimes_A M$,
\begin{align}
D\phi&=(D\circ Q_0)1=(Q_0\otimes \id)(S_m1)\label{D_q_0_etc}\\
&=(Q_0\otimes \id)(1\otimes m)=(Q_01)\otimes m\notag \\
&=\phi\otimes m.\notag 
\end{align}
Thus the map $Q_0$ is determined by a single vector $\phi\in X$ satisfying 
$\langle \phi, \phi \rangle=1$ and $D\phi=\phi \otimes m$. Following the classical case, 
we say that $(X, D, \phi)$ is a \emph{scaling function} for the filter $m$. 

\begin{cor}\label{def_of_MRA}
Suppose $M$ is a correspondence over a unital $C^*$-algebra $A$ such that $M$ is essential as a left $A$-module, $m$ 
is a filter in $M$, and $(X, D, 
\phi)$ is a scaling function for $m$. Then there is a Hilbert module isomorphism 
$R$ of $M_\infty$ onto $M_\infty \otimes_A M$ such that $R\circ \iota^{k+1}=\iota^k \otimes 
\id$, and there is an inner-product preserving map 
$Q$ of $M_\infty$ into $X$ such that $(Q\circ \iota^0)1=\phi$ and 
$D\circ Q=(Q\otimes \id)\circ R$. The family 
\[
\{V_k:=Q(\iota^k(M^{\otimes k})): k\geq 0\}
\] 
is a projective multi-resolution analysis for the submodule $Q(M_\infty)$ of $X$, and 
$V_0$ is the free rank-one $A$-module generated by $\phi.$ 
\end{cor}

\begin{proof} The formula $Q_0a=\phi \cdot a$ defines a module map $A\to X$. 
Calculations like the ones giving \eqref{D_q_0_etc} show that 
$D\circ Q_0=(Q_0\otimes \id)\circ S_m$, and   
Theorem~\ref{limit_from_general_V}(b) gives an inner-product preserving map 
$Q:M_\infty \to X$ with the required properties.
\end{proof}

\section{Multi-resolution analyses from transfer operators}\label{section_Exel}

Our first applications involve correspondences built from transfer 
operators for endomorphisms of $C^*$-algebras. Suppose $\alpha$ is an endomorphism of a unital $C^*$-algebra $A$. A positive linear map $L:A\to A$ is a \emph{transfer operator} for $\alpha$ if $L(a\alpha(b)) 
=L(a)b$ for $a,b\in A$. 

In \cite{exel}, Exel constructs a crossed product $A\times_{\alpha,L}\N$ using
a correspondence $M_L$ over~$A$. To construct $M_L$, he
endows the vector space $A_L:=A$ with the right action of~$A$ given by 
$m\cdot  a :=m\alpha(a)$ and the pre-inner product 
\begin{equation}
\langle m_1, m_2\rangle
=L(m_1^*m_2)\ \text{ for $m_1,m_2\in A_L$,}\label{def_of_inner_prod}
\end{equation}
and completes to get  a right Hilbert 
$A$-module $M_L$. The completing process includes modding out the vectors of length zero, and since $\|m\cdot 1-m\|^2$ is always zero, we have $m\cdot 1=m$ for every $m\in M_L$, so that $M_L$ is essential as a right $A$-module. 
The action of $A$ by left multiplication on $A_L$ extends to a left action of $A$ on $M_L$, which is implemented by a 
homomorphism $A\to \mathcal{L}(M_{L})$.  Exel's module $M_L$ is also essential 
as a left $A$-module, so the right and left module actions induce isomorphisms 
$M_L\otimes_A A\cong M_L \text{ and }A\otimes_A M_L \cong M_L$. 

We now present some simple lemmas which will help us work with the modules $M_L$. The first concerns the powers $L^k$, which are easily seen to be transfer operators for the powers $\alpha^k$ of the endomorphism $\alpha$. This lemma is essentially the same as Proposition~2.1 of \cite{Lar}, but the conventions there are a little different.

\begin{lemma}\label{M_Lk}
For each $k,l\in \N$, the map  $a\otimes b\mapsto a\alpha^k(b)$ of $A\odot A$ into $A$ induces an isomorphism of the  correspondence $M_{L^k}\otimes_A M_{L^l}$ onto $M_{L^{k+l}}$. 
\end{lemma}
 
\begin{proof}
It is easy to check that the map is a bimodule homomorphism which preserves the inner products,
 and hence extends to an injection of correspondences. To see that it  has dense range and is 
therefore surjective, note that for $a\in A_{L^{k+l}}$ we have $a=a\cdot 1=a\alpha^{k+l}(1)$ in $M_{L^{k+l}}$, and hence $a$ is the image of $a\otimes \alpha^l(1)\in M_{L^k}\otimes_A M_{L^l}$.
\end{proof}

\begin{lemma}\label{X=XL}
Suppose that $L$ is a transfer operator for $\alpha\in \End A$, and $X$ is a Hilbert $A$-module. 
\smallskip

\textnormal{(a)} The underlying vector space of $X$ becomes a 
pre-inner-product $A$-module with $x\cdot_L a:=x\cdot\alpha(a)$ and $\langle x,y\rangle_L:=L(\langle x,y\rangle)$. 
The completion is a Hilbert $A$-module, which we denote by $X_L$.

\smallskip

\textnormal{(b)} The map $\Phi$ of $X\odot A_L$ into $X$ defined in terms of the right action 
of $A$ on $X$ by $\Phi(x\otimes a)=x\cdot a$ induces an isomorphism of $X\otimes_A M_L$ onto $X_L$.
\smallskip

\textnormal{(c)} Suppose that $T:Y\to Z$ is an adjointable Hilbert 
$A$-module homomorphism. Then $T_L:=T\otimes \id :Y_L=Y\otimes_A M_L\to 
Z\otimes_A M_L
=Z_L$ is given by the formula $T_Ly=Ty$.
\end{lemma}

\begin{proof}
The claim of (a) follows from the defining properties of a transfer operator. For (b), 
just verify that $\Phi$ is inner-product preserving and has dense range. 
To prove (c), note that every $y\in Y_L$ has the form
$\Phi(y\otimes 1)$, and so $T_Ly=\Phi(T\otimes \id (y\otimes 1))=Ty$.
\end{proof}

Suppose that $C$ is a compact space, $\sigma:C\to C$ is a surjective local homeomorphism, and $\alpha:f\mapsto 
f\circ \sigma$ is the associated endomorphism of $C(C)$. As in \cite{ev}, the formula
\begin{equation}\label{exelsto}
L(f)(c)=\frac{1}{|\sigma^{-1}(c)|}\sum_{\sigma(d)=c}f(d)
\end{equation}
defines a transfer operator for $(C(C),\alpha)$. The choice of $|\sigma^{-1}(c)|^{-1}$ as normalising factor is not important: for any continuous function $w:C\to (0,\infty)$,
\[
L_w(f)(c)=w(c)\sum_{\sigma(d)=c}f(d)
\]
is also a transfer operator for $\alpha$. (The function $c\mapsto |\sigma^{-1}(c)|$ is locally constant and hence continuous.) We usually use \eqref{exelsto} because for this choice, filters in $M_L$ include the filters of wavelet theory (see Example~\ref{classicex} below). However, we need the extra generality in the next lemma because we want to apply it to the powers $L^k$ of $L$, which are transfer operators for $\alpha^k:f\mapsto f\circ \sigma^k$, but which need not have the obvious normalising factors $|\sigma^{-k}(c)|^{-1}$ (see Lemma~\ref{computeLk}).

\begin{lemma}\label{equivnorm}
For every $(C(C),\alpha,L_w)$ as above, and every Hilbert $C(C)$-module $X$, the function $\|\cdot\|_w$ on $X$ defined by
\begin{equation}\label{defLnorm}
\|x\|^2_w:=\|L_w(\langle x,x\rangle)\|=\sup_{c\in C}\Big(w(c)
\sum_{\sigma(d)=c}\langle x,x\rangle(d)\Big)
\end{equation}
is a norm which is equivalent to the given norm on $X$.
\end{lemma}

\begin{proof}
We first notice that $\|\cdot\|_w$ is the seminorm associated to the pre-inner product 
$\langle\cdot,\cdot\rangle_{L_w}$ used to define $X_{L_w}$, and is in particular a seminorm. Choose $\delta$, $M$ and $K$ such that $0<\delta\leq w(c)\leq M$ and $|\sigma^{-1}(c)|\leq K$ for every $c\in C$. Then we trivially have 
$\|x\|^2_w\leq MK\|x\|^2$. On the other hand, 
since $\langle x,x\rangle$ is a continuous non-negative function on a compact space, 
there exists $c\in C$ such that 
\[
\langle x,x\rangle(c)=\|\langle x,x\rangle\|_\infty =\|x\|^2,
\]
and then 
\begin{align*}
\|x\|_w^2&\geq L_w(\langle x,x\rangle)(\sigma(c))=w(\sigma(c))\textstyle{\sum_{\sigma(d)=\sigma(c)}}\langle x,x\rangle(d)\\
&\geq \delta\langle x,x\rangle(c)=\delta\|x\|^2.
\end{align*}
From this estimate we deduce, first, that $\|x\|_w=0$ implies $x=0$, so that $\|\cdot\|_w$ 
is a norm, and, second, that $\|\cdot\|_w$ is equivalent to the given norm.
\end{proof}

We now compute an explicit formula for $L^k$.

\begin{lemma}\label{computeLk}
The $k$th power of the transfer operator $L$ defined in \eqref{exelsto} is given by
\[
L^k(f)(c)=\sum_{\sigma^k(d)=c}\Big(\prod_{j=1}^k|\sigma^{-1}(\sigma^j(d))|^{-1}\Big)f(d).
\]
\end{lemma}

\begin{proof}
By induction on $k$. For $k=1$, the normalising factor is the one in \eqref{exelsto} because $\sigma(d)=c$. For the inductive step, we 
write
\[
L^{k+1}(f)(c)=\frac{1}{|\sigma^{-1}(c)|}\sum_{\sigma(d)=c}L^k(f)(d)=\sum_{\sigma(d)=c}|\sigma^{-1}(\sigma(d))|^{-1}L^k(f)(d)
\]
and expand $L^k(f)(d)$.
\end{proof}

Lemma~\ref{equivnorm} implies that, for the systems of the form $(C(C),\alpha,L)$, every 
Hilbert $C(C)$-module $X$ is already complete in the norm $\|\cdot\|_L$ used to define $X_L$, and $X_L$ has $X$ as its underlying vector 
space. Thus the isomorphism $D:X\to X\otimes_{C(C)} M_L$ which we use in Theorem~\ref{limit_from_general_V} to identify the direct limit $Y_\infty$ is in particular a linear isomorphism of $X$ onto $X$. 
(Though it is only a 
module homomorphism when we use the correct module action of ${C(C)}$ on $X_L$: we have~$D(x\cdot f)=(Dx)\cdot_L f=(Dx)\cdot\alpha (f)$.) This leads to more familiar-looking reformulations of Theorem~\ref{limit_from_general_V} and Corollary~\ref{def_of_MRA}.

\begin{prop}\label{dilationforCT}
Consider a system of the form $(C(C),\alpha,L)$, and Exel's correspondence $M_L$. Suppose $Y$ is a Hilbert $C(C)$-module and $T:Y\to Y_L$ is adjointable and inner-product preserving. Then $T$ is also adjointable and inner-product preserving as a map of $Y_{L^k}$ into $Y_{L^{k+1}}$. Let $(Y_\infty,\iota^k)=\varinjlim(Y_{L^k},T)$, and denote by $R:Y_\infty\to (Y_\infty)_L$ the map obtained from the isomorphism of Theorem~\ref{limit_from_general_V}(a) by identifying $Y_\infty\otimes_{C(C)} M_L$ with $(Y_\infty)_L$.

Now suppose that $X$ is a Hilbert $C(C)$-module, that $D:X\to X$ is a linear isomorphism of $X$ onto $X$ satisfying
\begin{equation}\label{propsdilv1}
D(x\cdot f)=(Dx)\cdot \alpha(f) \quad \text{and} 
\quad L(\langle Dx,Dy\rangle)=\langle x,y\rangle,
\end{equation}
and that $Q_0:Y\to X$ is inner-product preserving and satisfies $D\circ Q_0=Q_0\circ T$. Then there is an inner-product preserving map $Q$ of $Y_\infty$ into $X$ such that $Q\circ \iota^0=Q_0$ 
and $D\circ Q=Q\circ R$, and the submodules $V_k:=Q(\iota^k(Y_{L^k}))$ of $Q(Y_\infty)$ satisfy:
\begin{enumerate}
\item $V_0$ is a complemented $C(C)$-submodule of $Q(Y_\infty)$;
\smallskip
\item $V_k=D^{-k}(V_0)$ for $k\geq 0$;
\smallskip
\item $V_k$ is a complemented submodule of $V_{k+1}$ for $k\geq 0$;
\smallskip
\item $\bigcup_{k=0}^\infty V_k$ is dense in $Q(Y_\infty)$.
\end{enumerate}
\end{prop}

\begin{proof}
It is easy to check by writing $L^{k+1}$ as $L^k\circ L$ that $T:Y_{L^k}\to Y_{L^{k+1}}$ has the required properties.
The equations \eqref{propsdilv1} say that $D$ 
is a Hilbert-module isomorphism of $X$ onto $X_L$. The isomorphism of $X\otimes_{C(C)} 
M_L$ onto $X_L$ carries $x\otimes f$ to $x\cdot f$, and hence converts $Q_0\otimes\id$ into $Q_0:Y_L\to X_L$; thus the equation $D\circ Q_0=Q_0\circ T$ says that $D\circ Q_0=(Q_0\otimes \id)\circ T$. Now part (b) of Theorem~\ref{limit_from_general_V} gives an inner-product preserving map $Q:Y_\infty\to X$ which satisfies $Q\circ \iota^0=Q_0$ 
and $D\circ Q=(Q\otimes \id)\circ R$; when we identify $X\otimes_{C(C)} M_L$ with $X_L$, the second equation becomes $D\circ Q=Q\circ R$.

Properties (i), (iii) and (iv)  follow immediately from the corresponding properties of $Y_k:=\iota^k(Y_{L^k})$ in Theorem~\ref{limit_from_general_V}. For (ii), note that the equations $R\circ \iota^{k+1}=\iota^k\otimes \id$ say that $R$ maps $Y_{k+1}$ onto $Y_k\otimes_{C(C)}M_L$; since the identification of $Y_\infty\otimes_{C(C)} M_L$ with $(Y_\infty)_L$ takes $y\otimes f$ into $y\cdot f$ (for the original action of $C(C)$ on $Y_\infty$), it takes  $Y_k\otimes_{C(C)} M_L$ onto $Y_k$. Thus, viewed as a map from $Y_\infty$ to $(Y_\infty)_L$, $R$ carries $Y_{k+1}$ onto $Y_k$. The relation $D\circ Q=Q\circ R$ therefore implies that $D(V_{k+1})=V_k$, so that, at least as vector spaces, $V_k=D^{-k}(V_0)$.
\end{proof}

\begin{cor}\label{scalingforCT}
Consider a system of the form $(C(C),\alpha,L)$, let $m$ be a filter in Exel's correspondence $M_L$, and let $M_\infty=\varinjlim(M_L^{\otimes k},T_k)$ be the direct limit of the system $(M_L^{\otimes k},T_k)$ defined by \eqref{new_Sk}. 
Suppose that $X$ is a Hilbert $C(C)$-module, $\phi\in X$ satisfies $\langle \phi,\phi\rangle =1$, 
and $D:X\to X$ is a linear isomorphism of $X$ onto $X$ such that
\begin{equation}\label{propsdil}
D(x\cdot f)=(Dx)\cdot \alpha(f),\quad L(\langle Dx,Dy\rangle)=\langle x,y\rangle,\quad \text{and} 
\quad D\phi=\phi\cdot m.
\end{equation}
Then $(X,D,\phi)$ is a scaling function for $m$, and there is a Hilbert-module isomorphism $Q$ of $M_\infty$ into $X$ such that $(Q\circ \iota^0)1=\phi$ 
and $D\circ Q=Q\circ R$. Moreover, if $V_0$ denotes the submodule $Q(\iota^0(C(C)))$, then  
\[
\{V_k=Q(\iota^k(M_L^{\otimes k}))=D^{-k}(V_0) : k\geq 0\}
\]
is a projective multi-resolution analysis for $Q(M_\infty)$.
\end{cor}

\begin{proof}
As in Corollary~\ref{def_of_MRA}, we define $Q_0:C(C)\to X$ by $Q_0(f)=\phi\cdot f$, and $Q_0$ is inner-product preserving because $\langle \phi,\phi\rangle=1_{C(C)}$. The isometry $T=S_m:C(C)\to C(C)=C(C)_L$ is given by $T(f)=m\alpha(f)$. Thus for $f\in C(C)$ we have
\begin{align*}
D\circ Q_0(f)&=D(\phi\cdot f)=(D\phi)\cdot \alpha(f)=(\phi\cdot m)\cdot \alpha(f)\\&=\phi\cdot(m\alpha(f))=\phi\cdot(Tf)=Q_0\circ T(f),
\end{align*}
and the result follows from Proposition~\ref{dilationforCT}.
\end{proof}

\begin{example}\label{classicex}
Let $N\geq 2$ be an integer, take $A=C(\T)$ and $\alpha(f)(z)=f(z^N)$. Then 
\begin{equation}\label{classicTf}
L(f)(z):=\frac 1N \sum_{w^N=z}f(w)
\end{equation}
defines a transfer operator for $\alpha$. This fits the above structure with $C=\T$ and 
$\sigma(z)=z^N$, and the module $M_L$ is the one discussed in Example~\ref{keyexonT}. A filter $m\in M_L$ is a function $m\in C(\T)$ satisfying
\begin{equation}
1=\langle m, m\rangle(z)=\frac 1N \sum_{w^N=z}\vert m(w)\vert^2\text{ for all  $z\in \T$,} \label{m_norm_unit}
\end{equation}
so that if also $m(1)=N^{1/2}$ then $m$ is a low-pass filter in the sense of wavelet 
theory. Provided $m$ is low-pass and sufficiently smooth near $1$ (see, for example,  
\cite[Lemma 5.37]{fra} or \cite[Theorem 5.1.3]{Bra-Jor-book}), the infinite product 
\[
\phi(t):=\prod_{k=1}^\infty \big(N^{-1/2}m(\exp(2\pi iN^{-k}t))\big)
\]
converges uniformly on compact subsets to a 
function $\phi$ satisfying the scaling equation
\begin{equation}\label{phiscales}
N^{1/2}\phi(Nt)=m(e^{2\pi it})\phi(t)\ \text{ for $t\in \R$,}
\end{equation}
and
\begin{equation}\label{phifilter}
\sum_{k\in \Z} |\phi(t-k)|^2=1 \text{ for all $t\in \R$}.
\end{equation}

The properties of the scaling function $\phi$ are naturally expressed in terms of the Hilbert 
$C(\T)$-module $\Xi$ of \cite[\S1]{PR2} (see the discussion around our Equations~\eqref{right_action_Xi} and~\eqref{def_ip_Xi}). Because it is locally the uniform limit of continuous functions, the scaling function $\phi$ is continuous, 
and \eqref{phifilter} implies that $\phi$ is an element of $\Xi$ satisfying $\langle\phi,\phi\rangle=1_{C(\T)}$. With $D_N:\Xi\to\Xi$ defined by $(D_N\xi)(t)=N^{1/2}\phi(Nt)$, the scaling equation \eqref{phiscales} becomes 
\[
(D_N\phi)(t)=m(e^{2\pi it})\phi(t)=(\phi\cdot m)(t)\ \text{ for all $t\in \R$.}
\] 
An easy calculation shows that $D_N(\xi\cdot f)=(D_N\xi)\cdot\alpha(f)$, and
\begin{align*}
L(\langle D_N\xi,D_N\eta\rangle)(e^{2\pi it})
&=\frac{1}{N}\sum_{j=0}^{N-1}\langle D_N\xi,D_N\eta\rangle(e^{2\pi i(t-j)/N})\\
&=\frac{1}{N}\sum_{j=0}^{N-1}\sum_{k\in\Z}\overline{(D_N\xi)((t-j)/N-k)}
(D_N\eta)((t-j)/N-k)\\
&=\sum_{j=0}^{N-1}\sum_{k\in\Z}\overline{\xi(t-(j+kN))}
\eta(t-(j+kN))\\
&=\sum_{\ell\in\Z}\overline{\xi(t-\ell)}
\eta(t-\ell)\\
&=\langle \xi,\eta\rangle(e^{2\pi it}),
\end{align*}
so that $D_N$ is an isomorphism of $\Xi$ onto $\Xi_L$.

Thus $(\Xi,D_N,\phi)$ satisfies the hypotheses of Corollary~\ref{scalingforCT}, and there is a Hilbert-module isomorphism $Q$ of $M_\infty$ onto a closed $C(\T)$-submodule of $\Xi$. The  subspaces $V_k:=Q(\iota^k(M_L^{\otimes k}))$ satisfy:
\begin{enumerate}
\item[(1)] $V_0$ is the free rank-one $A$-module generated by $\phi=(Q\circ\iota^0)1$;
\smallskip
\item[(2)] $V_k=D_N^{-k}(V_0)$ for $k\geq 0$;
\smallskip
\item[(3)] $V_k\subset V_{k+1}$ for $k\geq 0$.
\end{enumerate}
These directly imply properties (1), (2) and (3) of \cite[Definition~4]{PR2}. Now define $V_k$ by (2) for $k<0$. Then, since $V_0$ is projective, Proposition~13 of \cite{PR2} implies that $\bigcap_{k=-\infty}^{0} V_k=\{0\}$, which is property (5) of \cite[Definition~4]{PR2}. Since the filter $m$ is 
low-pass, so that $m(1)=N^{1/2}$, the scaling function $\phi$ satisfies $\phi(0)=1$, and hence it 
follows from Proposition~14 of \cite{PR2} that $\bigcup_{k=0}^\infty V_k$ is dense in $\Xi$, which is property (4) of \cite[Definition~4]{PR2}. Thus the subspaces $\{V_k:k\in\Z\}$ form a projective multi-resolution analysis for $\Xi$ in the sense of \cite[Definition~4]{PR2}.
\end{example}

\section{Frames and orthonormal bases}\label{frames}

Suppose that $X$ is a Hilbert module over a unital $C^*$-algebra $A$. We say that a countable subset $\{x_j:j\in J\}$ in $X$ is an \emph{orthonormal basis}\footnote{In our classic Example~\ref{classicex}, an orthonormal basis will consist of finitely many functions $x_j:\T\to\C$, and an engineer would call such a basis ``a filter bank with perfect reconstruction''. That filter banks fit naturally into the setting of Exel's correspondences has also been noticed by other researchers, including Ionescu and Muhly \cite{im}.} for $X$ if $\{x_j\}$ generates $X$ 
and $\langle x_j,x_k\rangle =\delta_{j,k}1_A$. By modifying the standard Hilbert-space argument we obtain the \emph{reconstruction formula}
\[
x=\sum_j x_j\cdot 
\langle x_j, x\rangle \ \text{for every $x\in X$,}
\]
where we are asserting that the partial sums for the series converge in norm to $x$.  Frames are, loosely speaking, sets $\{x_j\}$ which are not orthonormal but still satisfy the reconstruction formula. Frames are particularly interesting in the context of Hilbert modules, which need not have an orthonormal basis, but nearly always have frames \cite{fra-lar, rt}.

A countable subset $\{x_j:j\in J \}$ 
in a Hilbert module $X$ over a unital $C^*$-algebra $A$ is a \emph{Parseval frame} for $X$ if it satisfies the \emph{frame identity} 
\begin{equation}\label{def_smf}
\langle x,x\rangle=\sum_{j\in J}\langle x,x_j\rangle \langle  x_j,x\rangle\ \text{ for every $x\in X$,}
\end{equation}
where we require that the sum on the right-hand side converges in norm in $A$. (These are called ``standard modules frames'' in \cite{PR2} and ``standard normalised tight frames'' in \cite{fra-lar}.) Equivalently, $\{x_j\}$ is a Parseval frame if and only if the reconstruction formula 
\[
x=\sum_{j\in J} x_j\cdot 
\langle x_j, x\rangle \ \text{holds for every $x\in X$,}
\]
where again one is asserting that the sum converges in norm in $X$ (this is proved in 
\cite{fra-lar} and generalised in \cite[Theorem~3.4]{rt}). Although one is in principle more interested in the reconstruction formula, the frame identity is often easier to manipulate than the reconstruction formula because it involves series whose terms are positive elements of the $C^*$-algebra $A$. 

Our definition of Parseval frame is at first sight a little incomplete, since we have not specified any order on the index set $J$. So we pause to prove the following reassuring lemma, which tells us exactly what we need to check to see that a given countable set is a frame.  We denote by $F(J)$ the set of finite subsets of $J$, directed by inclusion.

\begin{lemma}\label{verifyframe}
Suppose that $X$ is a Hilbert module over a $C^*$-algebra, $\{x_j:j\in J\}$ is a countable subset of $X$, and $x\in X$. For each finite subset $F\in F(J)$, we write
\[
s_F:=\sum_{j\in F}\langle x,x_j\rangle\langle x_j,x\rangle.
\]
Suppose that
\begin{enumerate}
\item[\textnormal{(a)}] $s_F\leq \langle x,x\rangle$ for every $F\in F(J)$, and
\smallskip

\item[\textnormal{(b)}] for every $\epsilon>0$ there exists $F\in F(J)$ such that $\|\langle x,x\rangle-s_F\|<\epsilon$.
\end{enumerate} 
Then $\|\langle x,x\rangle-s_F\|\to 0$ as $F$ runs through $(F(J),\leq)$.
\end{lemma}

Provided we order $J$ so that the sets $F_k:=\{j\in J:j\leq k\}$ are cofinal in $F(J)$, Lemma~\ref{verifyframe} implies that 
the partial sums $s_k:=s_{F_k}$ converge to $\langle x,x\rangle$ in norm. This applies in particular to the order given by any enumeration of $J$, and also to the product order when $J=K\times L$ is a product of two ordered sets with the same property. Notice that convergence with respect to any such order implies (a) and (b), and hence convergence in one such order implies convergence in all others. (This equivalence for enumerations was noted in \cite{fra-lar}.)

\begin{proof}[Proof of Lemma~\ref{verifyframe}]
Fix $\epsilon>0$, and choose $F$ as in (b). Then for any $G\in F(J)$ such that $F\subset G$, we can apply (a) to $s_G$, and
\[
s_F\leq s_F+\sum_{j\in G\backslash F}\langle x,x_j\rangle\langle x_j,x\rangle=
s_G\leq \langle x,x\rangle;
\]
subtracting gives
\[
0\leq \langle x,x\rangle-s_G
\leq \langle x,x\rangle-s_F.
\]
But $0\leq b\leq c$ implies $\|b\|\leq \|c\|$, and hence 
\[
F\subset G\Longrightarrow \|\langle x,x\rangle-s_G\|\leq  \|\langle x,x\rangle-s_F\|<\epsilon.
\]
Thus $s_G\to \langle x,x\rangle$ in the given ordering on $F(J)$, as claimed.
\end{proof}

\begin{prop}\label{frameMinfty}
Suppose that $M$ is a correspondence over a unital $C^*$-algebra $A$, that $M$ is essential as a left $A$-module, and that $\{m_j:j\geq 0\}$ is a Parseval frame for $M$ such that $\langle m_0,m_0\rangle =1$. Define $T_k:M^{\otimes k}\to M^{\otimes(k+1)}$ by $T_kn=m_0\otimes n$, and let $(M_\infty,\iota^k)=\varinjlim(M^{\otimes k},T_k)$. Then
\begin{equation}\label{frameforMinfty}
\{\iota^0(1)\}\cup\big\{\iota^k\big(\textstyle{\bigotimes_{l=1}^k}m_{j_l}\big):k\geq 1,\ j_l\geq 0\text{ and } j_1>0\big\}
\end{equation}
is a Parseval frame for $M_\infty$. If $\{m_j\}$ is an orthonormal basis for $M$, then \eqref{frameforMinfty} is an orthonormal basis for $M_\infty$.
\end{prop}

If $\{e_i\}$ and $\{f_j\}$ are orthonormal bases for Hilbert spaces $H$ and $K$, then the tensor products $\{e_i\otimes f_j\}$ form an orthonormal basis for the tensor product $H\otimes K$. In Proposition~\ref{frameMinfty}, though, we are dealing with tensor products which are balanced over actions of a possibly non-commutative algebra $A$. So the next lemma is not quite as obvious as it might seem.

\begin{lemma}\label{tensorframe}
Suppose that $X$ is a Hilbert $A$-module and $M$ is a correspondence over $A$. If $\{x_i:i\in I\}$ and $\{m_j:j\in J\}$ are  Parseval frames for $X$ and $M$ respectively, then 
\[
\{x_i\otimes m_j:i\in I\text{ and }j\in J\}
\]
is a Parseval frame for $X\otimes_A M$.
\end{lemma} 

\begin{proof}
First let $y=x\otimes m$ be an elementary tensor in $X\otimes_A M$; we aim to verify that $\{x_i\otimes m_j\}$ has properties (b) and (a) of Lemma~\ref{verifyframe}. By the frame identity for $\{x_i\}$, there is a finite subset $F$ of $I$ such that we have a norm approximation
\begin{align}\label{usexi}
\langle y,y\rangle&=\langle x\otimes m,x\otimes m\rangle=\langle\langle x,x\rangle\cdot m,m\rangle\\
&\sim\sum_{i\in F}\big\langle\langle x,x_i\rangle\langle x_i,x\rangle\cdot m,m\big\rangle\notag\\
&=\sum_{i\in F}\big\langle\langle x_i,x\rangle\cdot m,\langle x_i,x\rangle\cdot m\big\rangle.\notag
\end{align}
Now we apply the frame identity for $\{m_j\}$ to the elements $\langle x_i,x\rangle\cdot m$ to find a single finite subset $G$ of $J$ such that
\begin{align}\label{usemj}
\langle y,y\rangle&\sim\sum_{i\in F}\sum_{j\in G}\big\langle\langle x_i,x\rangle\cdot m,m_j\big\rangle\big\langle m_j,\langle x_i,x\rangle\cdot m\big\rangle\\
&=\sum_{i\in F}\sum_{j\in G}\langle y,x_i\otimes m_j\rangle\langle x_i\otimes m_j,y\rangle.\notag
\end{align}
Thus the set $\{x_i\otimes m_j\}$ satisfies (b). However, in the calculations \eqref{usexi} and \eqref{usemj}, we can replace the approximations $\sim$ by inequalities $\geq$, and deduce that $\{x_i\otimes m_j\}$ also satisfies (a) for this choice of $y$. Thus  Lemma~\ref{verifyframe} implies that $\{x_i\otimes m_j\}$ satisfies the frame identity for elementary tensors $y$.

The frame identity extends to elements of the algebraic tensor product $X\odot M$ by linearity, but to see that it extends to elements of the completion $X\otimes_A M$ seems to require some more work. 
We choose an enumeration $y_k=x_{i(k)}\otimes m_{j(k)}$ for the countable set $\{x_i\otimes m_j\}$. Then the frame identity implies that for each $y\in X\odot M$, the sequence $\{\langle y,y_k\rangle\}$ belongs to the Hilbert module $l^2(A)$, with 
\begin{equation}\label{frametrans}
\big\langle\{\langle y,y_k\rangle\},\{\langle y,y_k\rangle\}\big\rangle=\langle y,y\rangle.
\end{equation}
Thus the formula $\F(y)=\{\langle y,y_k\rangle\}$ defines a linear map of $X\odot M$ into $l^2(A)$, which by \eqref{frametrans} and the polarisation identity is inner-product preserving. Since $l^2(A)$ is a Hilbert module,  $\F$ extends to an inner-product preserving map $\overline{\F}$ of the completion $X\otimes_A M$ into $l^2(A)$. Since the maps $\overline{\F}$, $\{a_k\}\mapsto a_k=\langle\{a_k\},e_l\rangle$, and $y\mapsto \langle y,y_k\rangle$ are continuous, we have $\langle\overline{\F}(y),e_l\rangle=\langle y,y_l\rangle$ for all $y$ in the completion $X\otimes_A M$, and $\overline{\F}(y)=\{\langle y,y_k\rangle\}$ for all $y\in X\otimes_A M$. In particular, $\overline{\F}(y)=\{\langle y,y_k\rangle\}$ belongs to $l^2(A)$, so
the series
\[
\sum_{k=1}^\infty \langle y,y_k\rangle\langle y_k,y\rangle
\]
converges in norm with sum $\langle \overline{\F}(y),\overline{\F}(y)\rangle=\langle y,y\rangle$ for every $y\in X\otimes_A M$. Since this last assertion implies that (a) and (b) hold for every $y\in X\otimes_A M$, we deduce that $\{x_i\otimes m_j\}$ satisfies the frame identity in the strongest possible sense.
\end{proof}

\begin{remark}
We have given a detailed proof of Lemma~\ref{tensorframe} because we have found these convergence issues a little slippery. 
To see why we think such detours might be necessary, observe that the calculations \eqref{usexi} and \eqref{usemj} show 
that, if we start with sequential frames $\{x_i\}$ and $\{m_j\}$, then for every elementary tensor we have 
\begin{equation}\label{condconv}
\langle x\otimes m, x\otimes m\rangle=\lim_{k\to \infty}\Big(\lim_{n\to\infty}\sum_{i=1}^k\sum_{j=1}^n\langle x\otimes m,x_i\otimes m_j\rangle\langle x_i\otimes m_j,x\otimes m\rangle\Big).
\end{equation}
The asymmetry in the definition of the balanced tensor product means the calculation has to be done this way round, and we 
don't see any reason why we should expect to be able to reverse the order of the limits in \eqref{condconv}. 
\end{remark}

\begin{proof}[Proof of Proposition~\ref{frameMinfty}]
We begin by observing that plugging $x=m_0$ into the frame identity \eqref{def_smf} for the frame $\{m_j\}$ gives
\[
\langle m_0,m_0\rangle = \langle m_0,m_0\rangle +\sum_{j\geq 1}\langle m_0,m_j\rangle \langle m_j,m_0\rangle;
\]
since each $\langle m_0,m_j\rangle \langle m_j,m_0\rangle$ is a positive element of the $C^*$-algebra $A$, it follows that
\begin{equation}\label{othogmj}
\langle m_0,m_j\rangle=0\ \text{ for every $j>0$.}
\end{equation}

Next, we observe that the argument of \cite[Theorem~2]{PR2} applies to the direct sum decomposition of Proposition~\ref{moddecomp}, and deduce that it suffices to check that for each $k\geq 1$,
\[
\big\{\iota^k\big(\textstyle{\bigotimes_{l=1}^k}m_{j_l}\big):j_l\geq 0\text{ and }j_1>0\big\}
\]
is a Parseval frame for
\[
Z_{k-1}:=Y_k\ominus Y_{k-1}=\iota^k(M^{\otimes k}\ominus T_{k-1}(M^{\otimes(k-1)})).
\]
The orthogonality relation \eqref{othogmj} implies that each $\bigotimes_{l=1}^km_{j_l}$ with $j_1>0$ belongs to the complement of $T_{k-1}(M^{\otimes(k-1)})=\{m_0\otimes n:n\in M^{\otimes(k-1)}\}$. An induction argument using Lemma~\ref{tensorframe} shows that
\[
\big\{\textstyle{\bigotimes_{l=1}^k}m_{j_l}:j_l\geq 0\big\}
\] 
is a Parseval frame for $M^{\otimes k}$, and then the reconstruction formula for this frame implies that
\[
\big\{\iota^k\big(\textstyle{\bigotimes_{l=1}^k}m_{j_l}\big):j_l\geq 0 \text{ and }j_1>0\big\}
\]  
satisfies the reconstruction formula in $Z_{k-1}$, and hence is a frame for $Z_{k-1}$. 

We can verify by direct calculation that when $\{m_j\}$ is orthonormal, so is each $\{\bigotimes_{l=1}^k m_{j_l}:j_l\geq 0\}$ for each fixed $k$. (It is crucial in this calculation that $M$ is essential, so that $1\cdot m=m$ for all $m\in M$.) Thus the last assertion follows from the orthogonality of the summands $Z_k$.
\end{proof}

Now we suppose that $(X,D,\phi)$ is a scaling function for $m_0$, and aim to apply the isomorphism $Q$ of Corollary~\ref{def_of_MRA} to obtain a frame for the submodule $Q(M_\infty)$ of $X$. The isomorphism $Q$ satisfies $Q\circ\iota^k=D^{-k}\circ(Q_0\otimes\id_k)$, where $Q_0:A\to X$ is characterised by $Q_01=\phi$. To apply $Q_0\otimes\id_k$ to $\bigotimes_{l=1}^k m_{j_l}$, we need to recall that the isomorphism of $A\otimes_A M^{\otimes k}$ onto $M^{\otimes k}$ takes $1\otimes n$ to $1\cdot n=n$, and thus
\[
(Q_0\otimes\id_k)\big(\textstyle{\bigotimes_{l=1}^k m_{j_l}}\big)=
(Q_0\otimes\id_k)\big(1\otimes\big(\textstyle{\bigotimes_{l=1}^k m_{j_l}}\big)\big)=\phi\otimes
\big(\textstyle{\bigotimes_{l=1}^k m_{j_l}}\big).
\]
Thus Proposition~\ref{frameMinfty} gives:

\begin{cor}\label{standard_module_frame}
Suppose that $\{m_j:j\geq 0\}$ is a Parseval frame for $M$ with $\langle m_0,m_0\rangle =1$, and $(X,D,\phi)$ is a scaling function for $m_0$. Then
\begin{equation}\label{frameforQMinfty}
\{\phi\}\cup\big\{D^{-k}\big(\phi\otimes\big(\textstyle{\bigotimes_{l=1}^k}m_{j_l}\big)\big):k\geq 1,\ j_l\geq 0\text{ and } j_1>0\big\}
\end{equation}
is a Parseval frame for $Q(M_\infty)$. It is orthonormal if $\{m_j\}$ is.
\end{cor}

\begin{example}
Packer and Rieffel proved that if $N\geq 2$, $\alpha\in \End C(\T)$ is defined by $\alpha(f)(z)=f(z^N)$, and $L$ is given by \eqref{classicTf}, then the module $M_L$ admits an orthonormal basis. Indeed, their \cite[Proposition~1]{PR1} is much more general than this: it applies to the 
endomorphism of $C(\T^n)=C(\R^n/\Z^n)$ induced by a dilation matrix $A\in M_n(\Z)$ with $|\det A|\geq 2$. The key question considered in \cite{PR1}, however, asks which individual low-pass filters $m_0$ can be extended to  an orthonormal basis $\{m_j\}$ for $M_L$. They gave positive answers to this question in \cite[Theorem~2]{PR1}, which says in particular that any low-pass filter on $\T^n$ for $n\leq 4$ will extend to an orthonormal basis.

However, Packer and Rieffel also show in \cite[\S4]{PR1} that there exists a low-pass filter $m_0$ on $\T^5$ for a dilation matrix $A$ with $\det A=3$ for which it is not possible to find filters $m_1$ and $m_2$ such that $\{m_0,m_1,m_2\}$ is an orthonormal basis for the corresponding 
$M_L$. The point of their construction is that for this $m_0$, the complement $S_{m_0}(C(\T^5))^\perp=(m_0\cdot C(\T^5))^\perp$ is a projective $C(\T^5)$-module which is not free. However, since it is a direct summand of a free module $M_L$ of rank $3$, it has a Parseval frame $\{m_1,m_2,m_3\}$ with three elements. (For example, if $\{e_1,e_2,e_3\}$ is an orthonormal basis for $M_L$, we can take $m_j=(1-S_{m_0}S_{m_0}^*)e_j$.) Then $\{m_0,m_1,m_2,m_3\}$ is a Parseval frame for $M_L$ which has $\langle m_0,m_0\rangle =1$ but is not orthonormal. So there was some point in working out Corollary~\ref{standard_module_frame} for frames as well as orthonormal bases.  
\end{example}

We now apply Corollary~\ref{standard_module_frame} to Exel's modules $M_L$ for the systems $(A,\alpha,L)=(C(C),\alpha,L)$ discussed in  Section~\ref{section_Exel}. In this case the modules $X\otimes_A M_L^{\otimes k}$ can all be realised on the same 
underlying vector space: to see this, we first realise $M_L^{\otimes k}$ as $M_{L^k}$ using Lemma~\ref{M_Lk}, and then use Lemma~\ref{X=XL} to view $X\otimes_A M_L^{\otimes k}$ as $X_{L^k}$, which by Lemmas~\ref{equivnorm} and~\ref{computeLk} has underlying space $X$. The identification of $M_L^{\otimes k}$ with $M_{L^k}$ is defined on elementary tensors by
\begin{equation}\label{idotimes}
m_{j_1}\otimes m_{j_2}\otimes \cdots\otimes m_{j_k}\mapsto m_{j_1}\alpha(m_{j_2})\cdots\alpha^{k-1}(m_{j_k}),
\end{equation}
where the product on the right is the product in the algebra $C(C)$, which is the underlying space of each $M_{L^k}$. The map $x\otimes m\mapsto x\cdot m$ defined using the original action of $A$ on $X$ is then an isomorphism of $X\otimes_A M_{L^k}$ onto $X_{L^k}$. Thus Corollary~\ref{standard_module_frame} gives:

\begin{cor}\label{frameforMLinfty}
Consider the system $(C(C),\alpha,L)$ associated to a surjective local homeomorphism $\sigma$ of a  
compact space $C$, and suppose that $\{m_j:j\geq 0\}$ is a Parseval frame for $M_L$ with $\langle m_0,m_0\rangle=1$. Suppose that $X$ is a Hilbert $C(C)$-module, that $\phi\in X$ satisfies $\langle\phi,\phi\rangle=1$, and that $D:X\to X_L$ is an 
isomorphism of Hilbert $C(C)$-modules such that $D\phi=\phi\cdot m_0$. Then 
\begin{equation}\label{frameforQMLinfty}
\{\phi\}\cup\big\{D^{-k}\big(\phi\cdot\big(\textstyle{\prod_{l=1}^k}\alpha^{l-1}(m_{j_l})\big)\big):k\geq 1,\ j_l\geq 0\text{ and } j_1>0\big\}
\end{equation}
is a Parseval frame for the submodule $Q((M_L)_\infty)$ of $X$ described in Corollary~\ref{scalingforCT}. If $\{m_j\}$ is an orthonormal basis for $M_L$, then \eqref{frameforQMLinfty} is an orthonormal basis for $Q((M_L)_\infty)$.
\end{cor}

In the classical situation of Example~\ref{classicex}, \cite[Proposition~1]{PR1} says that the module $M_L$ is free, and hence we can apply Corollary~\ref{frameforMLinfty} to the scaling function $(\Xi,D,\phi)$, and thereby find module bases for $\Xi$. For $N=2$, we can do the calculations explicitly, and it is interesting to compare the resulting module basis with the analogous basis for $L^2(\R)$ obtained by applying Mallat's construction.
 
\begin{example}\label{basisforXi}
Let $N=2$. A filter $m_0$ in the module $M_L$ of Example~\ref{classicex} is a quadrature mirror filter: it satisfies
\begin{equation}\label{defqmf}
|m_0(z)|^2+|m_0(-z)|^2=2\ \text{ for all $z\in \T$.}
\end{equation}
(These are slightly different from the filters used in \cite{lr}, where we normalised so the sum in \eqref{defqmf} was $1$.) 
We assume that $m_0$ is low-pass, so that there exists a scaling function $\phi$ with respect to the usual 
dilation operator, and Corollary~\ref{scalingforCT} gives an isomorphism $Q$ of $(M_L)_\infty$ onto the module $\Xi$ of Packer and Rieffel~\cite{PR2} (see Example~\ref{classicex}). In this situation, we can write down specific functions $m_1\in M_L$ such that $\{m_0,m_1\}$ is an orthonormal basis for $M_L$: $m_1(z):=z\overline{m_0(-z)}$ is the usual choice. 

After $\phi$, the next element of our basis is the function
\[
\psi(t)=D^{-1}(\phi\cdot m_1)(t)=2^{-1/2}\phi(2^{-1}t)m_1(e^{\pi it})=2^{-1/2}e^{\pi it}\phi(2^{-1}t)\overline{m_0(-e^{\pi it})}.
\]
This is the same as the function denoted by $\psi$ in \cite{lr}, and is the Fourier transform of the usual dyadic mother wavelet. 
In $L^2(\R)$, the functions given by 
$(\psi\cdot e_k)(t)=e^{2\pi ikt}\psi(t)$ form an orthonormal basis for the space $\overline W_0:=\overline{Q(Z_0)}$ 
(which is the Fourier transform of the space usually denoted $W_0$ in the wavelet literature); here the single function $\psi$ is a basis for the Hilbert $C(\T)$-module $W_0$. 
At the next stage, the Hilbert module $W_1:=Q(Z_1)$ is free of rank $2$, with module basis
\[
\psi_l(t)=D^{-2}(\phi\cdot(m_1\alpha(m_l)))(t)=2^{-1}\phi(2^{-2}t)m_1(e^{\pi it/2})m_l(e^{\pi it})=(D^{-1}\psi)(t)m_l(e^{\pi it})
\]
for $l=0,1$. To get a Hilbert-space basis for $\overline W_1\subset L^2(\R)$ from this module basis, we need to include both
\begin{align*}
(\psi_0\cdot e_k)(t)&=e^{2\pi ikt}(D^{-1}\psi)(t)m_0(e^{\pi it})=D^{-1}(\psi\cdot e_{2k})(t)m_0(e^{\pi it})
\\
\intertext{and} 
(\psi_1\cdot e_k)(t)&=e^{2\pi ikt}(D^{-1}\psi)(t)m_1(e^{\pi it})=D^{-1}(\psi\cdot e_{2k+1})(t)\overline{m_0(-e^{\pi it})},
\end{align*}
and the resulting basis for $\overline W_1$ is slightly different from the usual Hilbert space basis $\{D^{-1}(\psi\cdot e_k)\}$ for $\overline W_1$. In general, 
as a $C(\T)$-submodule of $\Xi$, the space $W_k:=Q(Z_k)$ is free of rank $2^k$.
\end{example}

\section{The projective multi-resolution analyses of Packer and Rieffel}\label{projPR}

Here we apply our constructions to the example studied in Sections 4 and 5 of \cite{PR2}. 
We begin by fixing integers $c$ and $d$ greater than $1$, and consider the system $(C(\T^2),\alpha,L)$ associated to the local homeomorphism $\sigma:(w,z)\mapsto (w^c,z^d)$. (We could handle $c,d\in \Z\setminus\{0,\pm1\}$ 
at the expense of adding lots of absolute values.)

We next fix integers $a$ and $q$ with $q>0$, and consider the Hilbert $C(\T^2)$-module
\begin{equation}
Y(q,a):=\{\xi:\T\times \R\to \C: \xi(z, t-1)=z^a\xi(z,t)\},
\label{Y_qa_as_space}
\end{equation}
with action $(\xi\cdot f)(z,t)=\xi(z,t)f(z,e^{2\pi i qt})$ and inner product
\begin{equation}
\langle \xi, \eta\rangle(z,t)=\sum_{k=0}^{q-1}\overline{\xi\Big(z, \frac{t-k}q\Big)}
\eta\Big(z, \frac{t-k}q\Big).
\label{ip_on_Y_qa}
\end{equation}
We aim to build projective multi-resolution analyses starting from isometries on the 
Hilbert module $Y(q,a)$, aiming for those constructed in \cite[\S5]{PR2}. Because our 
constructions require us to use only modules over $C(\T^2)$, we have had to use 
slightly different realisations of the Hilbert modules used in \cite{PR2}. However, the function $\Phi(F)(z,t):=F(z,qt)$ is an isomorphism of the Hilbert 
$C(\T^2)$-module $X(q,a)$ described in \cite[Proposition~17]{PR2} onto $Y(q,a)$ compatible with the analogous isomorphism of $C(\R^2/(\Z\times q\Z))$ onto $C(\T^2)$.

We need an inner-product preserving adjointable map $S:Y(q,a)\to Y(q,a)_L$. Since 
Lemma~\ref{equivnorm} implies that   
$Y(q,a)_L$ has the same underlying space as $Y(q,a)$, an adjointable operator must 
in particular map $Y(q,a)$ to itself. For such a map to be a module homomorphism, it 
must send $\xi\in Y(q,a)$ to something involving the function $(z,t)\mapsto \xi(z^c,dt)$;
 this function belongs to $Y(q,cda)$, and to get back into $Y(q,a)$, Packer and Rieffel 
multiply it by 
an element of $Y(q,(1-cd)a)$. With our normalisation, $Y(q,(1-cd)a)$ has the same 
underlying space as $Y(1,(1-cd)a)$, and it is an inner product for this latter module 
which turns out to be relevant. We believe it is an advantage of our approach that the condition on $m$ which 
makes $S_m$ an isometry is simply expressed in terms of an inner product; Packer and Rieffel can only describe their 
condition (which appears below as \eqref{mtildePRfilter}) as ``closely related to one of the standard equations that a low-pass filter must satisfy'' \cite[page~459]{PR2}.

\begin{prop} Suppose that $m\in Y(1,(1-cd)a)_L$ satisfies $\langle m, m\rangle_L 
=1$. Then 
\begin{equation}
(S_m\xi)(z,t)=m(z,t)\xi(z^c,dt)\label{def_of_S_PR2}
\end{equation}
defines an inner-product preserving adjointable map $S_m:Y(q,a)\to Y(q,a)_L$.
\end{prop}

\begin{proof}
Let $\xi\in Y(q,a)$. Then $S_m\xi$ is certainly continuous, and
\begin{align*}   
(S_m\xi)(z, t-1)&=m(z,t-1)\xi(z^c, d(t-1))=z^{(1-cd)a}m(z,t)(z^c)^{ad}\xi(z^c, dt)\\
&=z^am(z,t)\xi(z^c, dt)=z^a(S_m\xi)(z,t),
\end{align*}
so $S_m\xi\in Y(q,a)$. To find an adjoint for $S_m$, we let $\xi,\eta\in Y(q,a)$ and
$(z,t)\in\T\times \R$, and compute 
\begin{align*}
\langle S_m\xi,&\eta\rangle_L(z,t)=L(\langle S_m\xi,\eta\rangle)(z,t)\\
&=\frac 1{cd}\sum_{w^c=z}\sum_{j=0}^{d-1}\sum_{k=0}^{q-1}
\overline{(S_m\xi)\Big(w,\frac 1q \Big(\frac {t-j}{d}-k\Big)\Big)}\eta\Big(w,\frac 1q \Big(\frac {t-j}{d}-k\Big)\Big)\\
&=\frac 1{cd}\sum_{w^c=z}\sum_{j=0}^{d-1}\sum_{k=0}^{q-1}
\overline{m\Big(w,\frac{t-j-kd}{qd}\Big)\xi\Big(z,\frac {t-j-kd}{q}\Big)}\eta\Big(w, \frac {t-j-kd}{qd}\Big).
\end{align*}
For the next step, we note that
\begin{align*} 
\{kd+j: 0\leq j\leq d-1,\; 
&0\leq k\leq q-1\}=\{n: 0\leq n\leq qd-1\}\\
&=\{ql+k: 0\leq l\leq d-1,\;0\leq k\leq 
q-1\}. 
\end{align*}
Thus
\begin{align*}
\langle S_m\xi,&\eta\rangle_L(z,t)\\
&=\frac 1{cd}\sum_{k=0}^{q-1}\sum_{w^c=z}\sum_{l=0}^{d-1}
\overline{m\Big(w, \frac{t-(ql+k)}{qd}\Big)\xi\Big(z, \frac {t-(ql+k)}{q}\Big)} 
\eta\Big(w,\frac{t-(ql+k)}{qd}\Big)\\
&=\frac 1{cd}\sum_{k=0}^{q-1} \overline{\xi\Big(z, \frac {t-k}q\Big)}\bigg(
\sum_{w^c=z}\sum_{l=0}^{d-1}
\overline{m\Big(w, \frac{t-(ql+k)}{qd}\Big)}\overline{z^{al}}
\eta\Big(w,\frac{t-(ql+k)}{qd}\Big)\bigg).
\end{align*}
Thus the function $T_m\eta$ defined by
\begin{equation}\label{defadjsm}
(T_m\eta)(z,s)=\frac 1{cd}\sum_{w^c=z}\sum_{l=0}^{d-1}\overline{m\Big(w,
\frac{s-l}d\Big)}\overline{z^{al}}\eta\Big(w,\frac{s-l}d\Big)
\end{equation}
satisfies $\langle S_m\xi,\eta\rangle_L=\langle \xi,T_m\eta\rangle$. We need to show that the formula \eqref{defadjsm} defines a function $T_m$ from $Y(q,a)=Y(q,a)_L$ to $Y(q,a)$, and then this equation says that $S_m$ is adjointable with adjoint $S_m^*=T_m$ (see \cite[\S2.2]{RW}).

The right-hand side of \eqref{defadjsm} gives a well-defined continuous function $T_m\eta:\T\times \R\to \C$, so we need to check that $(T_m\eta)(z,s-1)=z^a(T_m\eta)(z,s)$. We compute:
\begin{align*}
(T_m\eta)(z,s-1)&=\frac 1{cd}\sum_{w^c=z}\sum_{l=0}^{d-1}\overline{m\Big(w,
\frac{s-(l+1)}d\Big)}\overline{z^{al}}\eta\Big(w,\frac{s-(l+1)}d\Big)\\
&=\frac 1{cd}\sum_{w^c=z}\sum_{l=1}^{d}\overline{m\Big(w,
\frac{s-l}d\Big)}\overline{z^{a(l-1)}}\eta\Big(w,\frac{s-l}d\Big).
\end{align*}
The $l=d$ summand in the right-hand side is
\[
\overline{m\Big(w,\frac{s-d}d\Big)}\overline{z^{a(d-1)}}\eta\Big(w,\frac{s-d}d\Big) 
=\overline{w^{(1-cd)a} m\Big(w,
\frac{s}d\Big)}\overline{z^{a(d-1)}}w^a\eta\Big(w,\frac{s}{d}\Big),
\]
which is exactly what we'd get for $l=0$. Thus we can replace $\sum_{l=1}^d$ by $\sum_{l=0}^{d-1}$, factor out $z^a$, and deduce that $(T_m\eta)(z,s-1)=z^a(T_m\eta)(z,s)$. Thus $S_m$ is adjointable with $S_m^*=T_m$ given by \eqref{defadjsm}.

For $m, n\in Y(1, (1-cd)a)_L$, we compute
\begin{align*}
(S_m^*S_n\xi)(z,t)=L(\langle m,n\rangle_{Y(1,(1-cd)a)})(z,t)\xi(z,t)=\langle m, n\rangle_L(z,t)\xi(z,t),
\end{align*}
which shows in particular that $S_m^*S_m=1$ when $\langle m,m\rangle_L=1$.
\end{proof}

We can now fix a unit vector $m\in Y(1,(1-cd)a)_L$ and apply  
Theorem~\ref{limit_from_general_V} to $S_m:Y(q,a)\to Y(q,a)_L=Y(q,a)\otimes_{C(\T^2)} M_L$. 
This gives a projective multi-resolution analysis of a direct limit module $Y(q,a)_{\infty}$. 
Next we want to use Proposition~\ref{dilationforCT} to identify $Y(q,a)_\infty$ 
with the module $\Xi$ considered in \cite{PR2}, and thereby obtain a projective 
multi-resolution analysis of $\Xi$. This requires restrictions on $m$. 

It was proved in \cite[\S5]{PR2} that there is a function $\tilde{m}$ in 
$X(q, (1-cd)a)$, which satisfies
\begin{equation}\label{mtildePRfilter}
\sum_{w^c=z}\sum_{k=0}^{d-1}\Big|\tilde m\Big(w,t+\frac{kq}{d}\Big)\Big|^2=1\ \text{ for all $(z,t)\in \T\times\R$},
\end{equation}
and for which there exists $\tilde{\sigma}$ in $\Xi$ satisfying 
\begin{equation}\label{sigmatildeunit}
\sum_{m,n\in\Z} |\tilde\sigma(x+m,y+qn)|^2=1\ \text{for $(x,y)\in \R^2$} 
\end{equation}
and
\begin{equation}\label{sigmatilde}
\tilde{\sigma}(cx,dy)=\tilde{m}(e^{2\pi ix},y)\tilde{\sigma}(x,y)\ \text{ for $(x,y)\in \R^2$.}
\end{equation}
It will also be important for us that the function $\tilde{\sigma}$ in \cite{PR2} satisfies
\begin{equation}\label{simgatilde0=1}
\tilde\sigma(0,0)=1.
\end{equation}
Recalling that $\Phi(F)(z,t)=F(z,qt)$ defines an isomorphism $\Phi$ of $X(q,(1-cd)a)$ onto $Y(q,(1-cd)a)$, 
we take $m:=\sqrt{cd}\,\Phi(\tilde{m})$ in $Y(q,(1-cd)a)=Y(1,(1-cd)a)$,
and define $\sigma\in \Xi $ by $\sigma(x,y):=\tilde{\sigma}(x,qy)$. Then 
\begin{equation*}
\langle m,m\rangle_L(z,e^{2\pi it})
=cd\,L\big(\langle\Phi(\tilde m),\Phi(\tilde m)\rangle\big)(z,e^{2\pi it})=\sum_{w^c=z}\sum_{k=0}^{d-1} \Big|\tilde m\Big(w,\frac{qt-qk}{d}\Big)\Big|^2,
\end{equation*}
which is identically $1$ by \eqref{mtildePRfilter}. A calculation using \eqref{sigmatildeunit} shows 
that $\sigma$ satisfies $\langle \sigma,\sigma\rangle=1$ for the $C(\T^2)$-valued inner product on $\Xi$, 
and from \eqref{sigmatilde} we deduce that $\sigma$ satisfies the scaling equation
\begin{equation}\label{scaling_eq_on_Y}
\sqrt{cd}\,\sigma(cx,dy)=\sqrt{cd}\,\tilde{\sigma}(cx,qdy)=\sqrt{cd}\,\tilde{m}(e^{2\pi i x}, qy)\tilde{\sigma}(x,qy)=m(e^{2\pi i x},y)\sigma(x,y).
\end{equation}

Define  $D:\Xi\to \Xi$ by $(D\xi)(x,y)=\sqrt{cd}\,\xi(cx,dy)$. Then $D$ is certainly linear, and it is an isomorphism because 
we can write down an inverse. A straightforward calculation shows that $D(\xi\cdot f)=(D\xi)\cdot\alpha(f)$, and for $\xi,\eta\in\Xi$, we have
\begin{align*}
\langle D\xi, D\eta\rangle_L(e^{2\pi i x}, e^{2\pi i y})
&=\frac 1{cd}\sum_{k=0}^{c-1}\sum_{l=0}^{d-1}\sum_{m,n\in \Z}(\overline{D\xi}D\eta)\Big(\frac{x-k}{c} -m,\frac{y-l}{d}-n\Big)\\
&=\frac 1{cd}\sum_{k=0}^{c-1}\sum_{l=0}^{d-1}\sum_{m,n\in \Z}cd\,
(\overline{\xi}\eta)(x-k-cm,y-l-dn)\\
&=\sum_{p,q\in \Z}(\overline{\xi}\eta)(x-p,y-q)\\
&=\langle \xi, \eta\rangle (e^{2\pi i x}, e^{2\pi i y}).
\end{align*}
Thus $D$ satisfies the hypotheses \eqref{propsdilv1} of Proposition~\ref{dilationforCT}.
 
For $\xi\in Y(q,a)$, we define $Q_0\xi:\R^2\to \C$ by
\[
(Q_0\xi)(x,y)=\xi(e^{2\pi i x}, y/q)\sigma(x,y/q).
\]
We need to show that $Q_0\xi\in\Xi$ and that $Q_0$ is inner-product preserving from $Y(q,a)$ to $\Xi$. To see this, let $\xi,\eta\in Y(q,a)$. Then since $\langle\sigma,\sigma\rangle=1$, we have
\begin{align*}
\langle \xi,\eta\rangle(e^{2\pi ix},e^{2\pi iy})
&=\sum_{k=0}^{q-1}(\overline{\xi}\eta)\Big(e^{2\pi ix},\frac{y-k}{q}\Big)\\
&=\sum_{k=0}^{q-1}(\overline{\xi}\eta)\Big(e^{2\pi ix},\frac{y-k}{q}\Big)\Big(\sum_{m,n\in\Z}\Big|\sigma\Big(x-m,\frac{y-k}{q}-n\Big)\Big|^2\Big).
\end{align*}
The function $\overline{\xi}\eta$ satisfies 
$(\overline{\xi}\eta)(z,t-n)=(\overline{\xi}\eta)(z,t)$ for $n\in \Z$, and hence we can pull $(\overline{\xi}\eta)(e^{2\pi ix},\frac{y-k}{q})$ inside the second sum to get
\[
\langle \xi,\eta\rangle(e^{2\pi ix},e^{2\pi iy})
=\sum_{k=0}^{q-1}\sum_{m,n\in\Z}(\overline{\xi}\eta)\Big(e^{2\pi ix},\frac{y-k}{q}-n\Big)\Big|\sigma\Big(x-m,\frac{y-k}{q}-n\Big)\Big|^2.
\]
Next, we observe that $\frac{y-k}{q}-n=\frac{y-(k+nq)}{q}$, write $l=k+nq$ and deduce that
\begin{align}
\langle \xi,\eta\rangle(e^{2\pi ix},e^{2\pi iy})\label{proveQ0iso}
&=\sum_{m,l\in\Z}(\overline{\xi}\eta)\Big(e^{2\pi ix},\frac{y-l}{q}\Big)\Big|\sigma\Big(x-m,\frac{y-l}{q}\Big)\Big|^2\\
&=\sum_{m,l\in\Z}(\overline{Q_0\xi}Q_0\eta)(x-m,y-l).\notag
\end{align}
When $\eta=\xi$, the right-hand side of \eqref{proveQ0iso} is $\sum_{m,l\in\Z}|Q_0\xi(x-m,y-l)|^2$; since the left-hand side $\langle \xi,\xi\rangle$ is continuous on $\T^2$, Equation~\eqref{proveQ0iso} implies that $Q_0\xi\in \Xi$. Now Equation~\eqref{proveQ0iso} (for distinct $\xi$ and $\eta$) says that $Q_0$ is inner-product preserving from $Y(q,a)$ to $\Xi$, as required.

We next need to check that $D\circ Q_0=Q_0\circ S_m$. For $\xi\in Y(q,a)$, we compute using the scaling equation~\eqref{scaling_eq_on_Y}:
\begin{align*}
((D\circ Q_0)\xi)(x,y)
&=\sqrt{cd}\, (Q_0\xi)(cx,dy)=\sqrt{cd}\,\xi(e^{2\pi i cx}, dy/q)\sigma(cx, dy/q)\\
&= m(e^{2\pi i x},y/q)\xi(e^{2\pi i cx}, dy/q)\sigma(x,y/q)\\
&=(S_m\xi)(e^{2\pi i x},y/q)\sigma(x,y/q)=Q_0(S_m\xi)(x,y).
\end{align*}

We have now verified all the hypotheses of Proposition~\ref{dilationforCT}, and deduce that there is an inner-product preserving map $Q:Y(q,a)_\infty\to \Xi$ such that $V_k:=D^{-k}(Q_0(Y(q,a)))$ have the properties (i)--(iv) of Proposition~\ref{dilationforCT}. However, we also know from \eqref{simgatilde0=1} 
that $\sigma(0,0)=\tilde\sigma(0,0)=1$, and since it is easy to find elements $\xi$ of $Y(q,a)$ such that $\xi(1,0)\not=0$, $Q(Y_0)$ contains functions which do not vanish at $(0,0)$. Thus Proposition~14 of \cite{PR2} 
implies that $\bigcup_{k=0}^\infty V_k$ is dense in $\Xi$. Since the extra 
requirement $\bigcap_{k=-\infty}^0V_k=\{0\}$ in \cite{PR2} is automatic for 
projective multi-resolution analyses of $\Xi$ (by \cite[Proposition~13]{PR2}), 
we have recovered \cite[Theorem~6]{PR2}. 

We hope that our derivation helps to make it clear where some of the hypotheses in \cite{PR2} come from, and why they are necessary. We do not claim to have substantially simplified the arguments, since we have been content to rely on key analytic results from \cite{PR2}.

\section{Modules arising from directed graphs}\label{graphex}

In this section we discuss a new family of examples based on directed graphs. We consider the module $M_L$ associated to the backward shift on the one-sided infinite-path space of the graph, and aim to realise the direct limit $(M_L)_\infty$ as a module of functions on the two-sided infinite-path space. Part of our construction is very general: direct limits of modules associated to inverse systems can often be realised as modules of functions on the inverse limit.

\subsection{Realisations as modules of functions on inverse limits}
Suppose that $r_k:C_{k+1}\to C_k$ is an inverse system of compact spaces in which each $r_k$ 
is a surjective local homeomorphism. Set $w_0=1$, and suppose we have functions $w_k\in C(C_k)$ satisfying 
the consistency condition
\begin{equation}\label{consistweights}
\sum_{r_k(d)=c}|w_{k+1}(d)|^2=|w_k(c)|^2\ \text{ for every $k\geq 0$ and every $c\in C_k$.}
\end{equation}
For $k\geq 1$, we write $r^{(k)}:=r_0\circ r_1\circ\cdots \circ r_{k-1}:C_k\to C_0$. 
Then we can make each $C(C_k)$ into a Hilbert $C(C_0)$-module $X_k$ by defining 
$(x\cdot f)(c)=x(c)f(r^{(k)}(c))$ and
\[
\langle x,y\rangle_k(c)=\sum_{r^{(k)}(d)=c}\overline{x(d)}y(d)|w_k(d)|^2,
\]
and then completing (if necessary). Equation~\eqref{consistweights} implies that the functions $r_k^*:C(C_k)\to C(C_{k+1})$ defined  
by $r_k^*x=x\circ r_k$ satisfy
\[
\langle r_k^*x,r_k^*y\rangle_{k+1}=\langle x,y\rangle_k,
\]
and thus $(X_k,r_k^*)$ is a direct system of Hilbert $C(C_0)$-modules.

Our goal is to identify the direct limit $X_\infty:=\varinjlim X_k$ as a module of 
functions on the inverse limit $C_\infty:=\varprojlim C_k$, which is itself a compact 
space. To define the $C(C_0)$-valued inner product on $C(C_\infty)$, we need some 
measures whose existence will follow from the following standard lemma. It is proved in \cite{blpr}, for example.

\begin{lemma}\label{abstractmeas}
Suppose that $r_k:C_{k+1}\to C_k$ is an inverse system of compact spaces with each $r_k$ 
surjective, and $\mu_k$ is a family of measures on $C_k$ such that $\mu_0$ is a 
probability measure and 
\begin{equation}\label{consistency}
\int (x\circ r_k)\,d\mu_{k+1}=\int x\,d\mu_k\ \text{ for $x\in C(C_k)$.}
\end{equation} 
Let $C_\infty=\varprojlim (C_k,r_k)$, and denote the canonical map from $C_\infty$ to $C_k$ 
by $\pi_k$. Then there is a unique probability measure $\mu$ on $C_\infty$ such that 
\[
\int (x\circ \pi_k)\,d\mu=\int  x\,d\mu_k\ \text{ for $x\in C(C_k)$.}
\]
\end{lemma} 

We now fix $c\in C_0$, and define measures $\mu^c_k$ on the fibres $(r^{(k)})^{-1}(c)$ by
\begin{equation}\label{defmutk}
\int x \,d\mu^c_k:=\sum_{r^{(k)}(d)=c} x(d)|w_k(d)|^2.
\end{equation}
The consistency condition \eqref{consistweights} implies that the family $\{\mu^c_k:
k\geq 0\}$ satisfies \eqref{consistency}, and hence Lemma~\ref{abstractmeas} gives a probability measure $\mu^c$ on the fibre 
\[
\pi_0^{-1}(c)=\varprojlim((r^{(k)})^{-1}(c),r_k)
\]
such that
\begin{equation}\label{defmut}
\int (x\circ\pi_k)\,d\mu^c=\sum_{r^{(k)}(d)=c} x(d)|w_k(d)|^2\ \text{ for $x\in C(C_k)$.}
\end{equation}

\begin{prop}\label{geninvlim}
With $(x\cdot f)(d)=x(d)f(\pi_0(d))$ and 
\[
\langle x,y\rangle (c):=\int_{\pi_0^{-1}(c)} \overline{x(d)}y(d)\,d\mu^c(d),
\]
$C(C_\infty)$ is an inner-product module over $C(C_0)$. We denote the completion by 
$X_\infty$. Then the maps $\pi_k^*:C(C_k)\to C(C_\infty)$ extend to inner-product 
preserving homomorphisms of $X_k$ into $X_\infty$, and $(X_\infty,\pi_k^*)$ is a direct limit for $(X_k,r_k^*)$.
\end{prop}

\begin{proof}
It follows from \eqref{defmut} that $\langle x,y\rangle$ is continuous when $x$ and $y$ belong to $\pi_k^*(C(C_k))$. The Stone-Weierstrass theorem implies that $\bigcup_{k\geq 0}
\pi_k^*(C(C_k))$ is uniformly dense in $C(C_\infty)$, and the map $x\mapsto \langle x,y\rangle$ is uniformly continuous, so it follows that 
$\langle x,y\rangle$ is continuous for every $x,y\in C(C_\infty)$. The algebraic properties are easy to check, so we can indeed complete $C(C_\infty)$ to get a Hilbert module $X_\infty$. 

The formula \eqref{defmut} implies that $\pi_k^*$ is inner-product preserving. Since 
\[
\pi_{k+1}^*\circ r_k^*=(r_k\circ\pi_{k+1})^*=\pi_k^*,
\]
and since $\bigcup_{k\geq 0}\pi_k^*(C(C_k))$ is dense in $X_\infty$, the maps $\pi_k^*$ induce an isomorphism of $\varinjlim X_k$ onto $X_\infty$.
\end{proof}

\subsection{Systems associated to directed graphs} Let $E$ be a finite directed 
graph, consisting of a set $E^0$ of vertices, a set $E^1$ of edges, and maps $r,s:E^1\to E^0$ which identify the range and source of edges. We assume throughout that $E$ has no sources: every vertex receives at least one edge. In general, our conventions about directed graphs are those of \cite{cbms}. A path of length $k\geq 1$ is a sequence $\nu=\nu_1\dots \nu_k$ of edges such that $s(\nu_i)=r(\nu_{i+1})$ for all $i$, and $E^k$ denotes the set of paths of length $k$. We denote by $E^\infty$ the set of right-infinite paths $c=c_0c_1c_2\cdots$, which have range $r(c):=r(c_0)$ but no source. The  
space $E^\infty$ is a closed subset of the product space $\prod_{k=0}^\infty E^1$, and is therefore a compact Hausdorff space in the product topology; the cylinder sets $Z(\nu):=\{c\in E^\infty: c_i=\nu_{i+1} 
\text{ for }i\leq k-1\}$ associated to finite paths $\nu$ form a base of compact-open sets for the topology on $E^\infty$. For $c\in E^\infty$ and $\nu\in E^k$, we denote by $\nu c$ the right-infinite path $\nu_1\nu_2\cdots \nu_k c_0c_1\cdots$.

For the rest of this section, we consider the system $(C(E^\infty),\alpha,L)$ associated to the backward shift on $E^\infty$ defined by $\sigma(c_0c_1c_2\cdots)=c_1c_2c_3\cdots$. The transfer operator $L$ is given by
\begin{equation}\label{defgraphL}
L(f)(c)=\frac{1}{|\sigma^{-1}(c)|}\sum_{\sigma(d)=c}f(d)=\frac{1}{|s^{-1}(r(c))|}\sum_{s(e)=r(c)}f(ec).
\end{equation}
A filter in the corresponding module $M_L$ is determined by a weighting on the edges of $E$. More precisely, for each vertex $v\in E^0$ we choose a vector $(w(e):e\in s^{-1}(v))$ in $\C^{s^{-1}(v)}$ such that $\sum_{e\in s^{-1}(v)}|w(e)|^2=|s^{-1}(v)|$, and then 
\begin{equation}\label{defgraphm}
m=\sum_{e\in E^1}w(e)\chi_{Z(e)}
\end{equation}
is a filter in $M_L$. Our goal is to identify the associated direct limit $M_\infty:=(M_L)_\infty$ as a concrete module of functions on the space $E^{(-\infty,\infty)}$ of doubly infinite paths $c=\cdots c_{-2}c_{-1}c_0c_1c_2\cdots$ in $E$. 

Our first step is to realise $M_L^{\otimes k}$ as the module $M_{L^k}$ associated to the transfer operator $L^k$ for $\alpha^k$. The following formula for $L^k$ follows from Lemma~\ref{computeLk}:
\begin{equation}\label{FormforLk}
L^k(f)(c)=\sum_{\{\nu\in E^k\,:\,s(\nu)=r(c)\}}\Big(\prod_{j=1}^k|s^{-1}(s(\nu_j))|^{-1}\Big)f(\nu c),
\end{equation}
and the natural identification of $M_L^{\otimes k}$ with $M_{L^k}$ is described in \eqref{idotimes}. Since $E$ has no sources, the coefficients in \eqref{FormforLk} are all non-zero, and Lemma~\ref{equivnorm} implies that the modules $M_{L^k}$ all have the same underlying space $C(E^\infty)$. The isometry $S_m\otimes \id_k: M_L^{\otimes k}\to M_L^{\otimes (k+1)}$ is given by $(S_m\otimes \id_k)(n)=m\otimes n$ (see the discussion around Equation~\eqref{new_Sk}), and this goes into $M_{L^{k+1}}$ as $m\alpha(n)$. Thus, viewed as an isometry $T_k:M_{L^k}\to M_{L^{k+1}}$, $S_m\otimes \id_k$ is given by the formula
\[
T_k(f)(c)=m(c)f(\sigma(c))=\sum_{e\in E^1}w(e)\chi_{Z(e)}(c)f(\sigma(c))=w(c_0)f(c_1c_2\cdots).
\]
We will identify $\varinjlim(M_{L^k},T_k)$ by replacing $(M_{L^k},T_k)$ with an isomorphic direct system of the form discussed in Proposition~\ref{geninvlim}.

For $k\geq 0$, we introduce the path spaces 
\[
E^{[-k,\infty)}:=\{c=c_{-k}\cdots c_{-1}c_0c_1\cdots:s(c_i)=r(c_{i+1})\text{ for all $i$}\},
\]
and the homeomorphisms $\sigma_k:E^\infty\to E^{[-k,\infty)}$ defined by
\[
\sigma_k(c_0c_1c_2\cdots)=d_{-k}d_{-k+1}d_{-k+2}\cdots\ \text{ where $d_{-k+j}:=c_j$}.
\]
With $r_k:E^{[-k-1,\infty)}\to E^{[-k,\infty)}$ defined by $r_k(c_{-k-1}c_{-k}\cdots)=c_{-k}c_{-k+1}\cdots$, we have an inverse system $(E^{[-k,\infty)},r_k)$ of compact spaces. The inverse limit $\varprojlim (E^{[-k,\infty)},r_k)$ is the two-sided infinite path space $(E^{(-\infty,\infty)},\pi_k)$, with
\[
\pi_k(\cdots c_{-2}c_{-1}c_0c_1c_2\cdots)=c_{-k}c_{-k+1}\cdots c_{-1}c_0c_1c_2\cdots.
\]
We now  define weight functions $w_k\in C(E^{[-k,\infty)})$ by 
\begin{equation}\label{defweights}
w_k(c):=\prod_{j=-k}^{-1}w(c_j)|s^{-1}(s(c_j))|^{-1/2},
\end{equation}
and, recalling that $(w(e))_{s(e)=v}\in\C^{s^{-1}(v)}$ has norm $|s^{-1}(v)|^{1/2}$, verify that
\begin{align*}
\sum_{r_k(d)=c}|w_{k+1}(d)|^2
&=\sum_{s(e)=r(c)} |w(e)|^2|s^{-1}(s(e))|^{-1}\Big(\prod_{j=-k}^{-1}|w(c_j)|^2|s^{-1}(s(c_j))|^{-1}\Big)\\
&=\Big(\sum_{s(e)=r(c)}|w(e)|^2\Big)|s^{-1}(r(c))|^{-1}|w_k(c)|^2\\
&=|w_k(c)|^2.
\end{align*}
For this system, the functions $r^{(k)}$ are the projections of $E^{[-k,\infty)}$ on $E^\infty=E^{[0,\infty)}$, so 
\[
(r^{(k)})^{-1}(c)=\{\sigma_k(\nu c):\nu\in E^k\text{ and }s(\nu)=r(c)\}.
\]
Thus the Hilbert $C(E^\infty)$-module $X_k$ of the previous subsection has underlying space $C(E^{[-k,\infty)})$, module action
\[
(x\cdot f)(c_{-k}c_{-k+1}\cdots)=x(c_{-k}c_{-k+1}\cdots)f(c_0c_1\cdots),
\]
and inner product
\[
\langle x,y\rangle_k(c)=\sum_{\{\nu\in E^k\,:\,s(\nu)=r(c)\}}\Big(\prod_{j=1}^k|w(\nu_j)|^2|s^{-1}(s(\nu_j))|^{-1}\Big)\overline{x(\sigma_k(\nu c))}y(\sigma_k(\nu c)).
\]
Proposition~\ref{geninvlim} describes $\varinjlim(X_k,r_k^*)$ as a completion of $C(E^{( -\infty,\infty)})$.

To relate this direct limit to $\varinjlim(M_{L^k},T_k)$, we define $V_k:X_k\to M_{L^k}$ for $k\geq 1$ by
\[
(V_kx)(c)=\Big(\prod_{j=0}^{k-1} w(c_j)\Big)x(\sigma_k(c))\ \text{ for $c\in E^\infty$};
\]
calculations show that, provided the weights $w(e)$ are all non-zero, $V_k$ is a Hilbert-module isomorphism of $X_k$ onto $M_{L^k}$. Then for $k\geq 1$, $x\in C(E^{[-k,\infty)})=X_k$ and $c\in E^\infty$, we have $r_k\circ \sigma_{k+1}=\sigma_k\circ \sigma$, and
\begin{align*}
T_k\circ V_k(x)(c)&=w(c_0)(V_kx)(\sigma(c))\\
&=w(c_0)\big(\textstyle{\prod_{j=1}^{k}} w(c_j)\big)x(\sigma_k\circ\sigma(c))\\
&=\big(\textstyle{\prod_{j=0}^{k}} w(c_j)\big)x(r_k\circ\sigma_{k+1}(c))\\
&=V_{k+1}(r_k^*x)(c).
\end{align*} 
Thus we have a commutative diagram
\begin{equation}\label{commdiagpath}
\xygraph{
{C(E^\infty)}="w0":[rr]{X_1}="w1"^{r_0^*}:[rr]{X_2}="w2"^{r_1^*}:[rr]{\cdots}="w3"^{r_2^*}
"w0":[uu]{C(E^\infty)}="v0"^{{\id}}:[rr]{M_L}="v1"^{T_0}:[rr]{M_{L^2}}="v2"^{T_1}:[rr]{\cdots}="v3"^{T_2}
"w1":"v1"^{V_1}"w2":"v2"^{V_2}
}\end{equation}
In other words, the $V_k$ form an isomorphism of direct systems. We deduce that the two systems have isomorphic direct limits. 

So we expect the projective multi-resolution analyses of $M_\infty$ to give projective multi-resolution analyses for the module $X_{\infty}(E):=X_\infty$ of Proposition~\ref{geninvlim}. It remains to identify the dilation operator and a function $\phi$ 
in $X_{\infty}(E)$ which form a scaling function together with $X_\infty(E)$.

\begin{prop}
Suppose that $\{w(e):e\in E^1\}$ are non-zero complex numbers such that $\sum_{s(e)=v}|w(e)|^2=|s^{-1}(v)|$ for every $v\in E^0$, define $w_k:E^{[-k,\infty)}\to \C$ by \eqref{defweights}, and define $m$ by \eqref{defgraphm}. Then $m$ is a filter in Exel's correspondence $M_L$; let $M_\infty=\varinjlim(M_L^{\otimes k},T_k)$ be the direct limit of the system $(M_L^{\otimes k},T_k)$ defined by \eqref{new_Sk}. Let $\mu^c$ be the measure satisfying \eqref{defmut}, and let $X_\infty(E)$ be the Hilbert $C(E^\infty)$-module obtained by completing $C(E^{(-\infty,\infty)})$ in the inner product defined by the measures $\mu^c$, as in Proposition~\ref{geninvlim}.
Let $h$ denote the backward shift homeomorphism on $E^{(-\infty,\infty)}$, and define $D:C(E^{(-\infty,\infty)})\to C(E^{(-\infty,\infty)})$ by
\[
(Dx)(c)=m(\pi_0(c))x(h(c)).
\]
Then $D$ extends to a linear isomorphism of $X_\infty(E)$ onto itself, $(X_\infty(E), D, 1)$ is a scaling function for $m$, and there is an isomorphism $Q$ of $M_\infty$ onto $X_\infty(E)$ such that $V_k:=Q(\iota^k(M_{L^k}))$ is a projective multi-resolution analysis for $X_\infty(E)$.
\end{prop}

\begin{proof}
We verify the hypotheses \eqref{propsdil} of Corollary~\ref{scalingforCT}. We first let $x\in C(E^{(-\infty,\infty)})$, $f\in C(E^\infty)$, and compute, observing that $\pi_0\circ h=\sigma\circ \pi_0$:
\begin{align*}
D(x\cdot f)(c)&=m(\pi_0(c))(x\cdot f)(h(c))=m(\pi_0(c))x(h(c))f(\pi_0(h(c)))\\
&=(Dx)(c)f(\sigma(\pi_0(c)))=(Dx)(c)\alpha(f)(\pi_0(c))\\
&=((Dx)\cdot\alpha(f))(c).
\end{align*}
Next, we consider $x=x_1\circ\pi_k$ and $y=y_1\circ\pi_k$, and compute:
\begin{align*}
L(\langle Dx,Dy\rangle)(c)
&=\frac{1}{|s^{-1}(r(c))|}\sum_{s(e)=r(c)}\langle Dx,Dy\rangle(ec)\\
&=\sum_{s(e)=r(c)}|s^{-1}(s(e))|^{-1}\int_{\pi_0^{-1}(ec)}(\overline{x}y)(h(d))|m(ec)|^2\,d\mu^{ec}(d).
\end{align*}
Since $m(ec)=w(e)$ and we can write $x\circ h=x_1\circ\pi_k\circ h=x_2\circ\pi_{k-1}$, $y\circ h=y_2\circ \pi_{k-1}$, this is
\[
\sum_{s(e)=r(c)}|s^{-1}(s(e))|^{-1}\int_{(r^{(k-1)})^{-1}(ec)} (\overline{x_2}y_2)(d)|w(e)|^2\,d\mu^{ec}_{k-1}(d).
\]
Next, observe that $(r^{(k-1)})^{-1}(ec)=\{\sigma_{k-1}(\nu ec):\nu\in E^{k-1}, s(\nu)=r(e)\}$, recall from \eqref{defmutk} the definitions of the measures $\mu^{ec}_{k-1}$, notice that $x_2\circ \sigma_{k-1}=x_1\circ \sigma_k$, and continue:
\begin{align*}
L&(\langle Dx,Dy\rangle)(c)\\
&=\sum_{s(e)=r(c)}|s^{-1}(s(e))|^{-1}\!\!\!\sum_{\nu\in E^{k-1}\cap s^{-1}(r(e))}\!\!(\overline{x_1}y_1)(\sigma_{k}(\nu ec))|w(e)|^2|w_{k-1}(\sigma_{k-1}(\nu ec))|^2\\
&=\sum_{s(e)=r(c)}\ \sum_{\nu\in E^{k-1}\cap s^{-1}(r(e))}(\overline{x_1}y_1)(\sigma_{k}(\nu ec))|w_k(\sigma_k(\nu ec))|^2\\
&=\sum_{\lambda\in E^k\cap s^{-1}(r(c))} (\overline{x_1}y_1)(\sigma_k(\lambda c))|w_k(\sigma_k(\lambda c))|^2\\
&=\int_{\pi_0^{-1}(c)}\overline{x(d)}y(d)\,d\mu^{c}(d),
\end{align*}
which is just $\langle x,y\rangle(c)$. Thus $D$ is isometric from $C(E^{(-\infty,\infty)})\subset X_\infty(E)$ to $X_\infty(E)_L$, and extends to a linear isometry on $X_\infty(E)$; since $m(\pi_0(c))=w(c_0)$ is never $0$,  every function of the form $x_1\circ \pi_k$ is in the range of $D$, and hence $D$ is 
surjective. 

Since the element $\phi=1$ trivially satisfies $D\phi=\phi\cdot m$, Corollary~\ref{scalingforCT} says that $(X_\infty(E),D,\phi)$ is a scaling function for $m$, and that there is a Hilbert-module isomorphism $Q$ of $M_\infty$ into $X_\infty(E)$. Since $D^{-1}$ maps $\pi_{k-1}^*(X_{k-1})$ onto $\pi_k^*(X_{k})$, and $X_0$ is in the range of $Q$, the range of $Q$ contains every $\pi_k^*(X_k)$ and hence is dense in $X_\infty(E)$. Since the range of every isometry is closed, we deduce that $Q$ is surjective.
\end{proof}

Now we want to describe the structure of the modules $V_0$ and $W_k:=V_{k+1}\ominus V_k$ in the direct sum decomposition $X_\infty(E)=V_0\oplus\big(\bigoplus_{k=0}^\infty W_k\big)$ obtained from  Proposition~\ref{moddecomp}. Since $V_0$ is isomorphic to $A_A=C(E^\infty)_{C(E^\infty)}$, it is free of rank one. In general, though, $W_0$ will not be free or even projective. To see why, notice that we will not be able to expand our weight function $w:E^1\to \C$ to an orthonormal basis for $M_L$ unless the spaces $\C^{s^{-1}(v)}$ all have the same dimension --- that is, unless each vertex emits the same number of edges. Nevertheless, we can still describe the module $W_0$ in a very concrete way. We need a general lemma.

\begin{lemma}\label{Smpis}
Suppose that $X$ is a Hilbert $A$-module and $m$ is an element of $X$ such that $p:=\langle m,m\rangle$ is a projection in $A$. Then $S_m:a\mapsto m\cdot a$ is an inner-product preserving map of $pA$ onto a complemented submodule of $X$, and $S_mS_m^*$ is the orthogonal projection of $X$ on $S_m(pA)$.
\end{lemma}

\begin{proof}
It is easy to check that $S_m$ is inner-product preserving on $pA$:
\[
\langle S_m(pa),S_m(pb)\rangle=\langle m\cdot(pa),m\cdot(pb)\rangle=a^*p\langle m,m\rangle pb=a^*pb=\langle pa,pb\rangle.
\]
Next, we verify by direct computation that $\|m\cdot p-m\|^2=0$, and deduce that $m\cdot p=m$. 
This formula implies that the adjoint $S_m^*$ of $S_m$ in $\L(A,X)$, which is given by $S_m^*x=\langle m,x\rangle$, satisfies
\[
p(S_m^*x)=p\langle m,x\rangle=\langle m\cdot p,x\rangle=\langle m,x\rangle=S_m^*x,
\]
and hence has range in $pA$. So $S_m^*$ is also an adjoint for $S_m:pA\to X$. Thus Lemma~\ref{Hilbert_module_isometry} implies that $S_m(pA)$ is a complemented submodule of $X$, and that $S_mS_m^*$ is the projection onto $S_m(pA)$.
\end{proof}

We now return to the task of identifying the module $W_0$ in the module decomposition of $X_\infty(E)$. We set $N=\max\{|s^{-1}(v)|:v\in E^0\}$, and choose functions $w_n:E^1\to \C$ for $1\leq n\leq N$ such that for every $v\in E^0$ the vectors
\[
\big\{(|s^{-1}(v)|^{-1/2}w_n(e))_{e\in s^{-1}(v)}:1\leq n\leq |s^{-1}(v)|\big\}
\]
form an orthonormal basis for $\C^{s^{-1}(v)}$, such that $w_n(e)=0$ for $n>|s^{-1}(s(e))|$, and such that $w_1$ is the weight function $w$ we used to define the filter $m$ in \eqref{defgraphm}. Define
\[
m_n=\sum_{e\in E^1}w_n(e)\chi_{Z(e)},
\]
so that in particular $m_1=m$. Then
\begin{align*}
\langle m_n,m_n\rangle(c)&=\frac{1}{|s^{-1}(r(c))|}\sum_{s(e)=r(c)}|m_n(ec)|^2=\frac{1}{|s^{-1}(r(c))|}\sum_{s(e)=r(c)}|w_n(e)|^2\\
&=\begin{cases}1&\text{if $n\leq |s^{-1}(r(c))|$}\\
0&\text{if $n> |s^{-1}(r(c))|$.} 
\end{cases}
\end{align*}
Thus $\langle m_n,m_n\rangle$ is the characteristic function $p_n:=\chi_{\{c\,:\,n\leq|s^{-1}(r(c))|\}}$, which is a projection in $A=C(E^\infty)$, and Lemma~\ref{Smpis} says that the operators $S_{m_n}$ are isomorphisms of $p_nC(E^\infty)$ onto complemented submodules $S_{m_n}S_{m_n}^*(M_L)$ of $M_L$. 

We claim that  $M_L=\bigoplus_{n=1}^NS_{m_n}S_{m_n}^*(M_L)$. Since the operators $S_{m_n}S_{m_n}^*$ are projections in the $C^*$-algebra $\L(M_L)$, it suffices to check that $\sum_{n=1}^NS_{m_n}S_{m_n}^*=1$. For $f\in M_L$ and $c\in E^\infty$, we have
\begin{align*}
\sum_{n=1}^NS_{m_n}S_{m_n}^*f(c)&=\sum_{n=1}^N m_n(c)\langle m_n,f\rangle(\sigma(c))\\
&=\sum_{n=1}^N w_n(c_0)\Big(\frac{1}{|s^{-1}(r(c_1))|}\sum_{s(e)=r(c_1)}\overline{m_n(e\sigma(c))}f(e\sigma(c))\Big)\\
&=\Big(\sum_{s(e)=r(c_1)}\sum_{n=1}^N |s^{-1}(r(c_1))|^{-1}w_n(c_0)\overline{w_n(e)}\Big)f(e\sigma(c)).
\end{align*}
Since the vectors $(|s^{-1}(r(c_1))|^{-1/2}w_n(e))_{s(e)=r(c_1)}$ for $1\leq n\leq |s^{-1}(r(c_1))|$ are the columns of a unitary matrix, the rows are also pairwise orthogonal, and hence the inside sum is $1$ when $e=c_0$ and $0$ otherwise. Thus 
\[
\sum_{n=1}^NS_{m_n}S_{m_n}^*f(c)=f(c_0\sigma(c))=f(c),
\]
and we have proved the claim.

In conclusion, then, we find that the module $V_1$ has the form
\[
V_1\cong \bigoplus_{n=1}^NS_{m_n}S_{m_n}^*C(E^\infty)\cong\bigoplus_{n=1}^Np_nC(E^\infty),
\]
and the module $W_0$ in the module decomposition of $X_\infty(E)$ satisfies
\[
W_0\cong \bigoplus_{n=2}^Np_nC(E^\infty).
\]

\begin{remark}\label{rightname?}
If every vertex in $E$ emits the same number of edges, then $p_n=1$ for every $n$, and the module $W_0$ is free. In general, though, $\{c:n\leq|s^{-1}(r(c))|\}$ will be a proper subset of $E^\infty$, $p_n$ will not be the identity of $C(E^\infty)$, and the modules in our ``projective multi-resolution analysis'' will not be projective. We have stuck with the name because we wanted to emphasise the connections with the work of Packer and Rieffel. However, our multi-resolution analysis for $X_\infty(E)$ seems to more closely resemble the generalised multi-resolution analyses of Baggett, Medina and Merrill \cite{bmm}, with the function $c\mapsto |s^{-1}(r(c))|$ playing the role of their multiplicity function. 
\end{remark}

\end{document}